\documentclass[11pt,leqno]{amsart}
\usepackage{enumerate,amsthm,amsmath,amssymb,comment,psfrag, float,dsfont}
\usepackage{xspace,tikz}
\usepackage{mathrsfs,mathtools,wasysym}
\usepackage{xcolor}
\usepackage{theoremref}
\usepackage{fancybox}
\usepackage{subfigure}
\usepackage[norelsize]{algorithm2e}
\usepackage{bm}
\usepackage{enumerate, verbatim}
\usepackage{stmaryrd}
\usepackage[top=1.5cm, bottom=1.5cm, outer=4cm, inner=3cm, heightrounded, marginparwidth=2.0cm, marginparsep=0.3cm]{geometry}
\input amssym.def
\numberwithin{equation}{section}

\begin{document}

\addtolength{\subfigcapskip}{-0.1cm}

\newtheorem{lemma}{Lemma}
\newtheorem{theorem}[lemma]{Theorem}
\newtheorem{remark}[lemma]{Remark}
\newtheorem{corollary}[lemma]{Corollary}
\newtheorem{notation}[lemma]{Notation}
\newtheorem{definition}{Definition}
\newtheorem{assumption}{Assumption}
\newtheorem{subsec} {} [section]
\newtheorem{example}{Example}
\vspace*{2.0cm}

\newcommand{\bb}{\boldsymbol{b}}

\newcommand{\nn}{\boldsymbol{n}}
\newcommand{\bu}{\boldsymbol{u}}
\newcommand{\bv}{\boldsymbol{v}}
\newcommand{\bw}{\boldsymbol{w}}
\newcommand{\bxi}{\boldsymbol{\xi}}
\newcommand{\ff}{\boldsymbol{f}}
\newcommand{\bx}{\boldsymbol{x}}

\newcommand{\btau}{\boldsymbol{\tau}}
\newcommand{\bvarphi}{\boldsymbol{\varphi}}

\newcommand{\bV}{\boldsymbol{V}}
\newcommand\calQ{\mathscr{Q}}
\newcommand{\bX}{\boldsymbol{X}}

\newcommand\calM{\mathscr{M}}
\newcommand\calF{\mathscr{F}}
\newcommand\calP{\mathscr{P}}
\newcommand\calC{\mathscr{C}}
\newcommand\calT{\mathscr{T}}
\newcommand\calL{\mathscr{L}}
\newcommand\calI{\mathscr{I}}
\newcommand\calB{\mathscr{B}}
\newcommand\calG{\mathscr{G}}

\newcommand\bR{\mathbb{R}}
\newcommand\bN{\mathbb{N}}
\newcommand\PP{\mathbb{P}}

\newcommand\sumMTH{\sum_{M\in\calT_H^{}}}

\newcommand\weakcv{ \rightharpoonup}

\newcommand{\dive}{{\ensuremath\mathop{\mathrm{div}\,}}}

\newcommand{\hypothesis}{{\ensuremath\mathop{\mathrm{Assumption~(A1)}}}}

\newcommand{\diam}{{\ensuremath\mathop{\mathrm{diam}\,}}}

\newcommand*\circled[1]{\tikz[baseline=(char.base)]{
            \node[shape=circle,draw,inner sep=2pt] (char) {#1};}}

\title[A stabilised finite element method for power-law fluids ]{Analysis of a stabilised finite element method for power-law fluids}


\author[G.R. Barrenechea]{Gabriel R. Barrenechea}
\address[G.R.B.]{Department of Mathematics and Statistics, University of Strathclyde, 26 Richmond Street, Glasgow G1 1XH, Scotland. {\tt gabriel.barrenechea@strath.ac.uk} }

\author[E. S\"uli]{Endre S\"uli}
\address[E.S.]{Mathematical Institute, Andrew Wiles Building, University of Oxford, Woodstock Road, Oxford OX2 6GG, UK.
{\tt endre.suli@maths.ox.ac.uk }}

\begin{abstract}
A low-order finite element method is constructed and analysed for an incompressible non-Newtonian flow problem with power-law rheology. 
The method is based on a continuous piecewise linear approximation of the velocity field and piecewise constant approximation of the pressure. Stabilisation, in the form
of pressure jumps, is added to the formulation to compensate for the failure of the inf-sup condition, and using an appropriate lifting of the
pressure jumps a divergence-free approximation to the velocity field is built and included in the discretisation of the convection term. This construction allows us to prove the convergence of the resulting finite element method for the entire range $r>\frac{2
d}{d+2}$ of the power-law index $r$ for which weak solutions to the model are known to exist in $d$ space dimensions, $d \in \{2,3\}$.
\end{abstract}


\maketitle

\begin{center}
\textit{Dedicated to Ron DeVore on the occasion of his 80th birthday\\~\\}
\end{center}

\section{Introduction}

The construction and mathematical analysis of finite element approximations of models of non-Newtonian fluids has been a subject of active research in recent years. Some of the most general results in this direction concern the convergence of mixed finite element approximations of models of incompressible fluids with implicit constitutive laws relating the Cauchy stress tensor to the symmetric velocity gradient (cf. \cite{DKS13}, \cite{ST20} and \cite{FGS}). Motivated by the groundbreaking contributions of  Cohen, Dahmen and DeVore \cite{CDD1,CDD2} and  Binev, Dahmen and DeVore \cite{BDD} concerning the convergence of adaptive algorithms for linear elliptic problems, progress, albeit much more limited in both scope and extent, has also been made on the analysis of adaptive finite element approximations of implicitly constituted non-Newtonian fluid flow models (cf. \cite{KS}).

Upon decomposing the Cauchy stress tensor into its traceless part, called the deviatoric stress tensor or shear-stress tensor, and its diagonal part, called the mean normal stress, models of incompressible fluids typically involve the velocity of the fluid, $\bu$,  its pressure, $p$, and the shear-stress tensor, $\boldsymbol{\mathcal{S}}$. For
Newtonian fluids the shear-stress tensor is a scalar multiple of the symmetric velocity gradient. The finite element approximation of Newtonian fluids is therefore 
usually performed in the velocity-pressure formulation.
For non-Newtonian fluids on the other hand the situation is more involved, because the shear-stress tensor exhibits nonlinear dependence as a function of the symmetric velocity gradient, and the functional relationship between the shear-stress tensor and the symmetric velocity gradient may even be completely implicit and multi-valued. For power-law fluids, such as the ones considered in this work, the shear-stress tensor exhibits power-law type growth as a function of the symmetric velocity gradient, the simplest instance of which results in an $r$-Laplace type operator in the balance of linear momentum equation, with a power-law exponent $r \in (1,\infty)$; for $r=2$, corresponding to a Newtonian fluid, the operator is linear, the Laplace operator. 
From a mathematical point of view, in the presence of a convection term in the balance of linear momentum equation in the model, the lower the value of $r$ the more difficult  the problem is to analyse. The existence of solutions for small values of 
$r$ was first proved  in \cite{FMS03}, where an Acerbi--Fusco type Lipschitz truncation was used in conjunction with Minty's method from monotone operator theory; thus, weak solutions were shown to exist
for $r>\frac{2d}{d+2}$ in $d\geq 2$ space dimensions. 

Finite element approximations of problems with power-law rheology have been extensively studied, including stabilised (or variational-multiscale) methods (cf. \cite{CC14,ACCCB18}, for example) and local discontinuous Galerkin methods (see, \cite{KRT14}, for example). The relevant literature is vast and it is beyond the scope of this work to provide an exhaustive survey of the various contributions; the interested reader may wish to consult \cite{OP02}, for example.  Concerning implicitly-constituted models, in the recent papers \cite{DKS13,ST20} the 
convergence of generic inf-sup stable velocity/pressure-based mixed finite element methods was proved for $r>\frac{2d}{d+1}$, while convergence for the full range, $r>\frac{2d}{d+2}$, was shown to be achievable only in the case the finite element methods where the velocity space consists of pointwise divergence-free functions. The reason for this dichotomy is that in the case of velocity approximations that are discretely divergence-free only, as is the case in generic inf-sup stable mixed finite element methods, the finite element approximation 
$\nabla \cdot(\bu_h \otimes \bu_h)$ of the convection term $\nabla \cdot(\bu \otimes \bu)$ does not vanish when tested with $\bu_h$, and it needs to be skew-symmetrized (cf. \cite{Tem77}) for this to happen. While in the case of the Navier--Stokes equations (corresponding  to $r=2$) membership of the velocity field to the natural function space for weak solutions, $W^{1,2}_0(\Omega)^d$, ensures that the convection term and its skew-symmetric modification can be bounded by the same expression using H\"older's inequality, this is not the case for the power-law model under consideration here for entire range $r>\frac{2d}{d+2}$ for which weak solutions to the problem are known to exist.
In fact, in the case of non-Newtonian power-law models the natural function space for the velocity field is $W^{1,r}_0(\Omega)^d$, and while the original convection term can be bounded in terms of the $W^{1,r}(\Omega)^d$ norm for all $r>\frac{2d}{d+2}$, for the skew-symmetric modification of the convection term, whose use is essential so as to be able to derive an energy inequality for discretely divergence-free velocity fields, this can only be achieved for the limited range $r > \frac{2d}{d+1}$. This was precisely the bottleneck encountered in \cite{DKS13} for discretely divergence-free velocity approximations, resulting in the reduction of the range of $r$ from the maximal range $r>\frac{2d}{d+2}$ for which weak solutions are known to exist, to $r>\frac{2d}{d+1}.$


The advantage of pointwise divergence-free finite element methods over discretely divergence-free finite element methods is therefore that, besides the physical consistency they provide, there is no need to rewrite the convection term in a skew-symmetric form. 
The topic of divergence-free finite element spaces has been treated extensively in the literature,
most commonly presenting pairs of spaces for which the divergence of the velocity space is a subspace of, or equal to, the pressure space. For example, the early Scott--Vogelius element \cite{SV84} (analysed recently in \cite{GS18})  uses $W^{1,2}(\Omega)^d$-conforming piecewise polynomials of degree $k$ for the velocity, while discontinuous piecewise polynomials of degree $k-1$ are used for the pressure. 
The stability of this pair requires either special meshes, or a high-enough degree $k$ (for example,  $k\ge 4$ is needed in \cite{GS18}, and $k=1$ is  only allowed in very special cases such as those described in Remark~\ref{Rem:low-order}). Another possibility is to relax the continuity requirements and consider a discontinuous Galerkin method, as was done, for example, in \cite{CKS07}, 
or to relax only the tangential continuity of the approximate velocity on faces of elements while still preserving its continuity in the direction of the normal to faces of elements, thus using $H(\dive\!;\Omega)$-conforming methods, as was the case in \cite{SL18}, for example. In this latter case the viscous term (defined as the divergence of the shear-stress $\mathcal{S})$ needs to be modified, for stability reasons, by adding terms controlling the jumps and averages of the velocity into the formulation, with, 
obviously undesirable, extra complications if the viscous term in the balance of linear momentum equation has a more complex structure, as is the case for the power-law model considered herein.

The recent works \cite{ABN18,ABN21} offer a way of preserving the advantages of a pointwise divergence-free approximation to the velocity field while working with the, computationally simplest, lowest-order  $W^{1,2}(\Omega)^d/L^2(\Omega)$-conforming velocity/pressure pair, namely,  $\mathbb{P}_1^d\times\mathbb{P}_0^{\rm disc}$. The key idea in those works can be summarised as follows: the discrete continuity equation 
contains a stabilising term based on the jumps of the discrete pressure. As the jumps of the pressure are constant along element faces, there exists a unique Raviart--Thomas field
such that its normal component is equal to the jumps. This field can be built at no extra computational cost, and then the continuity equation can be rewritten as a standard continuity equation,
but for a modified velocity field, which is now solenoidal. The finite element method then involves replacing the original discrete velocity field $\bu_h$ with the new, now solenoidal, modified velocity field
in the convection term. This facilitates the proofs of stability and convergence of the resulting finite element method without the need to rewrite the convection term in a skew-symmetric
form.  Our aim here is to apply this idea to a problem in non-Newtonian fluid mechanics. As a first step in this direction, we have chosen an explicit constitutive law  with power-law rheology.
 Even though this is the simplest
 constitutive law, it has been shown experimentally to faithfully reproduce many situations of physical interest (see the discussion in \cite{GRC21}, and the experimental results in, e.g., \cite{HE+11});
we therefore believe that it is a representative model for exemplifying the applicability of the proposed method in a mathematically nontrivial and physically relevant setting. 
Since the convection term does not need to be rewritten in a skew-symmetric form, the resulting method can now be proved to be stable and convergent to a weak solution for the whole range $r>\frac{2d}{d+2}$ of the power-law index for which weak solutions to the model are known to exist. In addition,
the sequence of numerical approximations is shown to converge strongly, and this strong convergence result is, to the best of our knowledge, a new contribution even in the, very special, Newtonian case ($r=2$). 

The rest of the manuscript is organised as follows. A section on preliminaries, containing the necessary notational conventions, basic definitions and results, the finite element spaces,
the lifting operator, the definition of the stabilising form, and properties of the discrete Lipschitz truncation method that we use, are presented following this Introduction. An important ingredient enabling the use of the discrete Lipschitz truncation technique is a discrete inf-sup condition that is given in the Appendix.
The finite element method is presented in Section~\ref{Sec:FEM}, where
we also show a uniform boundedness result for the sequence of approximations. Based on this and results pertaining to the discrete Lipschitz truncation, in Section~\ref{Sec:Convergence} the convergence
of the discrete solution to a weak solution of the model problem is proved using a compactness argument. Finally, some conclusions are drawn and potential future extensions are indicated.


\section{Preliminaries}\label{Sec:Prelim}

\subsection{Notation and the problem of interest}  

We use standard notation for Sobolev spaces. In particular, for $D\subset\bR^d$, $d=2,3$ and $s\in [1,+\infty]$, we denote  by $W^{k,s}_0(D)$ the closure of $C_0^\infty(D)$ with respect to the $W^{k,s}(D)$  norm, and
by $L^s_0(D)$ the space of functions in $L^s(D)$ with zero mean value. 
The norm in $L^s(D)$ is denoted by $\|\cdot\|_{0,s,D}^{}$; when $s=2$ we shall use the simpler notation $\|\cdot\|_{0,D}^{}$, and
the inner product in $L^2(D)$ will be denoted by $(\cdot,\cdot)_D^{}$.
For $k\ge 0$, the norm (seminorm) in $W^{k,s}(D)$
is denoted by $\|\cdot\|_{k,s,D}^{}$ ($|\cdot |_{k,s,D}^{}$). Moreover, the space $W^{-1,s'}(D)$ is the dual of $W^{1,s}_0(D)$ 
with duality pairing denoted by $\langle\cdot,\cdot\rangle_D^{}$.
We also denote by $W^s(\dive\!;D)$ the space of functions in $L^s(D)^d$ whose distributional
divergence belongs to $L^s(D)$, and by $W_0^{s}(\dive\!;D)$ the set of elements in $W^s(\dive\!;D)$ whose normal
trace on $\partial D$ is zero. 
In the above inner products and norms we do not make a distinction between scalar- and vector- or tensor-valued
functions.


Let $\Omega\subset \mathbb{R}^d, d=2,3,$ be an open, bounded, polyhedral domain with a Lipschitz boundary.
In this work we treat the problem with power-law rheology: given $r\in (1,\infty)$ and a right-hand side $\ff\in W^{-1,r'}(\Omega)^d$, 
find the velocity $\bu$, the pressure $p$, and the shear-stress tensor $\boldsymbol{\mathcal{S}}$ satisfying
\begin{equation}\label{Strong}
\left\{
\begin{array}{rcll}
-\dive \boldsymbol{\mathcal{S}} +\dive (\bu\otimes\bu)+\nabla p 
& = &
\ff & \text{in }\;\Omega,\\
\dive \bu & = & 0 &\text{in }\;\Omega,\\
\bu & = & \boldsymbol{0} &\text{on }\partial\Omega.
\end{array}
\right.
\end{equation}
There are many possible choices for the constitutive law, linking $\boldsymbol{\mathcal{S}}$ and the velocity $\bu$. In this work we have chosen the power-law description where
$\boldsymbol{\mathcal{S}} = \eta\, |\nabla\bu|^{r-2}\nabla\bu$, where $\eta>0$ is a reference viscosity. In order to simplify matters we will suppose that $\eta=1$,
but we should keep in mind that, to maintain physical consistency this reference value should be kept. Similarly, in physically realistic models the gradient of the velocity is usually
replaced by the symmetric velocity gradient $\varepsilon(\bu):=\frac{1}{2}(\nabla\bu +\nabla\bu^t)$. The results obtained in this paper can be extended, with minor modifications based on Korn's inequality, to that case as well, so for the sake of simplicity of the exposition we shall proceed with the constitutive relation 
$\boldsymbol{\mathcal{S}} = \eta\, |\nabla\bu|^{r-2}\nabla\bu$ (with $\eta = 1$) instead of 
$\boldsymbol{\mathcal{S}} = \eta\, |\varepsilon(\bu)|^{r-2}\varepsilon(\bu)$.

In order to state the weak formulation of \eqref{Strong} we need to present a few additional ingredients associated with the exponent in the constitutive law relating
$\boldsymbol{\mathcal{S}}$ and $\bu$. For $r\in (1,\infty)$, let $r'$ be its conjugate given by the relation $\frac{1}{r}+ \frac{1}{r'}=1$, and
let us define the critical exponent $\tilde{r}$ as follows:
\begin{equation}\label{rtilde-definition}
\tilde{r}:=\min\left\{r',\frac{r^\star}{2}\right\},\quad\textrm{where}\quad r^\star:=\left\{
\begin{array}{cl} \infty & \textrm{if}\; r\ge d, \\ \dfrac{dr}{d-r} & \textrm{otherwise}. \end{array}\right.
\end{equation}

\begin{figure}[H]
\centering
\includegraphics[width=7.0cm,height=5.0cm]{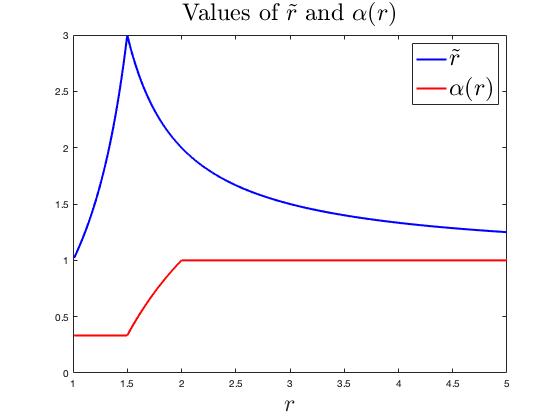}
\includegraphics[width=7.0cm,height=5.0cm]{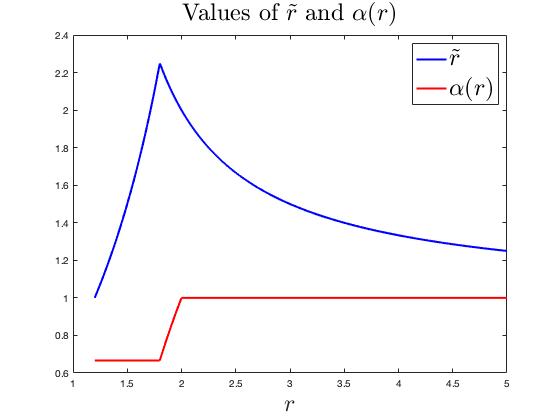}
\caption{Values of $\tilde{r}$ (defined in \eqref{rtilde-definition}) and $\alpha(r)$ (defined in \eqref{tauF-non-Newtonian}) for the cases $d=2$ (left) and $d=3$ (right).}
\label{Fig1}
\end{figure}

\begin{remark} 
With the definition \eqref{rtilde-definition} of $\tilde{r}$, the space $W^{1,r}(\Omega)$ is continuously embedded
in $L^{r^\star}(\Omega)$ if $r<d$ and in $L^s(\Omega)$, for every $s<\infty$, if $r\ge d$ (see, e.g., \cite[Corollary~9.14]{Brezis}). Then, 
in particular, $W^{1,r}(\Omega)$ is continuously embedded in $L^{2\tilde{r}}(\Omega)$ and there exists a $C>0$ such that
\begin{equation}\label{Sobolev-emb}
\|v\|_{0,2\tilde{r},\Omega}^{}\le C\,\|v\|_{1,r,\Omega}^{}\qquad \forall\, v\in W^{1,r}(\Omega).
\end{equation}

Moreover, the value of $\tilde{r}$ exhibits two different regimes, as can be seen in Figure~\ref{Fig1}, where its range of values is depicted. We will distinguish between $\tilde{r}\le 2$ and $\tilde{r} > 2$. The latter case occurs
for $r\in \big( \frac{4d}{d+4},2\big)$ and the maximum value of $\tilde{r}$  is attained when $r'=\frac{r^\star}{2}$, at which point we have the following values:
\begin{equation}\label{rtil-max}
r=\frac{3d}{d+2} = \left\{ \begin{array}{ll} \frac{3}{2} & \textrm{if}\; d=2,\\~ \\ \frac{9}{5} &  \textrm{if}\; d=3, \end{array} \right. \qquad\textrm{and}\qquad 
\tilde{r}_{\rm max}^{} = \frac{3d}{2d-2} = \left\{ \begin{array} {ll} 3 & \textrm{if}\; d=2,\\ \frac{9}{4} &  \textrm{if}\; d=3. \end{array} \right.
\end{equation}
\end{remark}

With this choice of stress tensor $\mathcal{S}$, the weak formulation of \eqref{Strong} is as follows:  find $\bu\in W^{1,r}_0(\Omega)^d$ and 
$ p\in L^{\tilde{r}}_0(\Omega)$ such that
\begin{alignat}{2}
(|\nabla\bu|^{r-2}\nabla\bu,\nabla\bv)_\Omega^{}-(\bu\otimes\bu,\nabla\bv)_\Omega^{}-(p,\dive\bv)_\Omega^{}&=
\langle \boldsymbol{f},\bv\rangle_\Omega^{}\qquad &&\forall\,\bv\in W^{1,\tilde{r}'}_0(\Omega)^d, \label{r-lap-1}\\
(q,\dive\bu)_\Omega^{} &= 0\qquad &&\forall q\in L^{r'}_0(\Omega). \label{r-lap-2}
\end{alignat}

\begin{remark}\label{rem1-1}
In order for the variational formulation \eqref{r-lap-1}, \eqref{r-lap-2} to be meaningful it is necessary that $\bu \otimes \bu \in L^{\tilde r}(\Omega)^{d \times d}$ with $\tilde r>1$, which necessitates that $r> \frac{2d}{d+2}$, and under this condition the existence of a solution to \eqref{r-lap-1}, \eqref{r-lap-2} has been proved (see \cite{DMS08}).
Thus, for the rest of this work we will assume that  $r> \frac{2d}{d+2}$.

Another fundamental ingredient in the proof of existence of solutions to \eqref{r-lap-1}, \eqref{r-lap-2} is the following inf-sup condition
(for a proof, see \cite{Galdi-VolI}): for $s,s'\in (1,+\infty)$ satisfying
$\frac{1}{s}+\frac{1}{s'}=1$, there exists a constant $\beta_s^{}>0$ such that
\begin{equation}\label{inf-sup-s}
\sup_{\bv\in W^{1,s'}_0(\Omega)^d\setminus\{\boldsymbol{0}\}}\frac{( q,\dive\bv)_\Omega^{}}{|\bv|_{1,s',\Omega}^{}}\ge \beta_s^{}\|q\|_{0,s,\Omega}^{}\qquad\forall\, q\in L^{s}_0(\Omega).
\end{equation}
\end{remark}


\subsection{Finite element spaces and preliminary results}\label{FEspaces}

Let $\{\calT_h^{}\}_{h>0}^{}$ be a shape-regular family of triangulations of $\overline{\Omega}$ consisting of closed
 simplices $K$ of diameter $h_K^{}\le h:=\max\{ h_K^{}:K\in\calT_h^{}\}$. To avoid technical difficulties we will suppose that the
family of triangulations is quasi-uniform. For reasons that will become apparent later, in the proof of convergence of the finite element method 
we will distinguish
between the cases $r\ge \frac{3d}{d+2}$ and $r\in \big( \frac{2d}{d+2}, \frac{3d}{d+2}\big)$. To cover the latter case (and for that purpose only) we need
to make the  following assumption
on the mesh:
\vspace{.25cm}

\noindent \textbf{Assumption (A1).} The triangulation $\calT_h^{}$ is the result of performing one (for
$d=2$), or two (for $d=3$), red refinement(s) of a, coarser, shape-regular triangulation $\calT_{H}^{}$. 
\vspace{.25cm}

We will denote the (closed) elements contained in $\calT_H^{}$   (referred to, in some instances, as {\sl macro}-elements) by $M$.

\begin{remark}\label{Rem:low-order}
$\hypothesis$ raises the question whether the space $\mathbb{P}_1^d\times\mathbb{P}_0^{}$ itself is stable on carefully constructed meshes. Some result are known in this direction. For example, 
in two space dimensions, this pair is inf-sup stable on Powell--Sabin meshes \cite{Zhang2008}. 
In the recent work  \cite{guzman_neilan_2}   local inf-sup stability is proved for this element in 
  barycentrically-refined meshes (also known as the
Alfeld split \cite{lai_schumaker_2007}, and a Hsieh--Clough--Tocher triangulation, see the references quoted in  \cite[P.~461]{Zhang2008}).
This is then used to build enriched elements that are divergence-free (although the velocity space contains quadratic face bubbles with constant divergence).
 For three-dimensional meshes, for the Alfeld split the lowest order inf-sup stable pair is $\mathbb{P}_4^3\times\mathbb{P}_3^{\rm disc}$ (cf. \cite{Zhang2005}), 
while for the Powell--Sabin split the lowest order inf-sup stable pair is $\mathbb{P}_2^3\times\mathbb{P}_1^{\rm disc}$  \cite{Zhang}.
 However, for the case considered in this paper, that is, taking the $\mathbb{P}_1^d\times\mathbb{P}_0^{}$ pair on general shape-regular meshes, stabilisation is a necessity.
In addition, it is important to note  that the papers cited above concern the Newtonian case only  and are mostly focused on the Stokes equations. The analysis of some of those
alternatives in the case of non-Newtonian flow models treated in the present work has not been carried out so far, and it will constitute a topic of future research.
\end{remark}

\begin{remark} 

\noindent (i) By letting $H:=\max\{ \diam(M):M\in\calT_H^{}\}$, clearly, $H\le Ch$, where $C$ does not depend on $h$. In fact, $C=2$ for $d=2$ and
$C=4$ for $d=3$.

\noindent (ii) Under $\hypothesis$, for every $\tilde{F}$, a facet of $M\in\calT_H^{}$, there exists at least one node  of $\calT_h^{}$ that belongs to the interior of $\tilde{F}$. In fact, this last remark is the main reason why $\hypothesis$ has been made on the meshes. In particular,  $\calT_h^{}$ could also result from first making a barycentric refinement of each facet of $\calT_H^{}$ and then building a conforming triangulation of $\overline\Omega$. For ease of exposition we shall simply adopt $\hypothesis$ in what follows.
\end{remark}

In the triangulation $\calT_h^{}$ we shall use the following notation:

\begin{itemize}
\item $\calF_h^{}$ : the set of all facets (edges in
$2D$ and faces in $3D$) of the triangulation $\calT_h^{}$, with diameter  $h_F^{}:=\diam(F)$.  
The set of internal facets is denoted by $\calF_I^{}$ and those on the boundary of $\Omega$ are denoted
by $\calF_{\partial}^{}$, so $\calF_h^{}=\calF_I^{}\cup\calF_\partial^{}$;

\item for every $M\in\calT_H^{}$ we denote by $\calF_I^{}(M)$ the set of facets of $\calT_h^{}$ whose interior
lies in the interior of $M$;

\item  for $F\in\calF_h^{}$ and $K\in \calT_h^{}$ we define the neighbourhoods
\begin{equation}\label{neighbourhoods}
\omega_F^{}:=\{K\in\calT_h^{}:F\in \calF_K^{}\},\quad  \quad \omega_K^{}:=\{K'\in\calT_h^{}:K\cap K'\not=\emptyset\};
\end{equation}

\item for each facet $F\in\calF_I^{}$ and every piecewise regular function $q$, we denote by 
$\llbracket q\rrbracket_F^{}$ the jump of $q$ across $F$;

\item for $\ell\ge 0$ we denote by $\PP_\ell^{}(K)$
the space of polynomials defined on $K$ of total degree smaller than, or equal to, $\ell$, and introduce
the following finite element spaces:
\begin{align}
\bV_h^{} &:=\{ \bv_h^{}\in C^0(\overline{\Omega})^d \,:\, \bv_h^{}|_K^{}\in\PP_1^{}(K)^d\;,\;\forall\, K\in\calT_h^{}\;,\;\bv_h^{}|_{\partial\Omega}^{}=\boldsymbol{0}\},\label{Vh}\\
\calQ_h^{} &:= \{ q_h^{}\in L^1_0(\Omega)\,:\, q_h^{}|_K^{}\in\PP_0^{}(K)\;,\;\forall\, K\in\calT_h^{}\},\label{Qh}\\
\calQ_H^{} &:= \{ q_H^{}\in L^1_0(\Omega)\,:\, q_H^{}|_M^{}\in\PP_0^{}(M)\;,\;\forall\, M\in\calT_M^{}\}.\label{QH}
\end{align}

\end{itemize}


Using the finite element spaces defined in \eqref{Vh}--\eqref{QH}, we denote by $S_h^{}:W^{1,r}_0(\Omega)^d\to \bV_h^{}$ the Scott--Zhang interpolation operator and by 
$\Pi_h^{}:L^1_0(\Omega)\to \calQ_h^{}$, $\Pi_H^{}:L^1_0(\Omega)\to \calQ_H^{}$  the projections defined by (see, e.g.,
\cite{EG21-I}):
\begin{align}
\Pi_h^{}q|_K^{} = \frac{(q,1)_K^{}}{|K|}\qquad\forall\; K\in \calT_h^{}, \label{def-Pih}\\
\Pi_H^{}q|_M^{} = \frac{(q,1)_M^{}}{|M|}\qquad\forall\; M\in \calT_H^{}. \label{def-PiH}
\end{align}
These operators satisfy (\cite{EG21-I}):
\begin{align}
\lim_{h\to 0} S_h^{}\bv =\bv\qquad\textrm{strongly in}\; W^{1,s}_0(\Omega)^d\,\quad\textrm{for all}\; \bv\in W^{1,s}_0(\Omega)^d\quad\textrm{and all}\; s 
\in [1,\infty), \label{conver-SZ}\\
\lim_{H\to 0} \Pi_H^{}q=\lim_{h\to 0} \Pi_h^{}q = q\qquad\textrm{strongly in}\; L^{s}_0(\Omega)^d\,\quad \textrm{for all}\; q\in L^{s}_0(\Omega)\quad\textrm{and all}\; s\in [1,\infty). \label{conver-Pih}
\end{align}

The following result, whose proof can be carried out using the techniques presented in \cite[Lemma~2.23]{EG21-I},  will be fundamental in the derivation (and analysis) of the proposed finite element method: for 
every $s\in (1,\infty)$ there exists a constant $C_s^{}>0$, independent of $h$, such that
\begin{equation}\label{average-bound}
\|q_h^{}-\Pi_H^{}(q_h^{})\|_{0,s,M}^{}\le C_s^{}\left\{\sum_{F\in\calF_I^{}(M)}h_F^{}\|\llbracket q_h^{}\rrbracket\|_{0,s,F}^s \right\}^{\frac{1}{s}},
\end{equation}
for all $M\in \calT_H^{}$, all $q_h^{}\in \calQ_h^{}$,  and all $h>0$.

We now recall three inequalities that will be useful in what follows. Let $s\in (1,\infty)$, $F\in \calF_h^{}$ and $K\in\omega_F^{}$. 
The following local trace inequality is a corollary of  the multiplicative trace inequality 
proved in \cite[Lemma~12.15]{EG21-I}:
\begin{equation}
\|v\|_{0,s,F}^{} \le C\, (h_F^{-\frac{1}{s}}\|v\|_{0,s,K}^{}+h_F^{1-\frac{1}{s}}\|\nabla v\|_{0,s,K}^{})\,.  \label{discrete-trace}
\end{equation}
In addition, we recall the following local inverse inequality (see, e.g., \cite[Lemma~12.1]{EG21-I}): for all $m,\ell\in\mathbb{N}, m\le \ell$ and all $p,q\in[1,+\infty]$, there
exists a constant $C$, independent of $h$, such that
\begin{equation}\label{inverse}
\|q\|_{\ell,p,K}^{}\le Ch_K^{m-\ell+d\left(\frac{1}{p}-\frac{1}{q}\right)}\|q\|_{m,q,K}^{},
\end{equation} 
for every polynomial function $q$ defined on $K$. A global version of this inequality can also be derived using the quasi-uniformity of the mesh family.
Finally, for $1<s\le \tilde{s}\le\infty$, a set of indices $\mathcal{I}$, and any vector $(x_i^{})_{i\in\mathcal{I}}^{}\in \ell^{\tilde{s}}(\mathcal{I})$,  the following inequality holds
(see \cite[Proposition~3.4(a)]{DFGHS03}  for its proof):
\begin{equation}\label{bound-lp}
\left\{\sum_{i\in \mathcal{I}}x_i^{\tilde{s}}\right\}^{\frac{1}{\tilde{s}}}\le \left\{\sum_{i\in \mathcal{I}}x_i^{s}\right\}^{\frac{1}{s}}.
\end{equation}

Finally, we note that under $\hypothesis$ the spaces $\bV_h^{}$ and $\calQ_H^{}$ satisfy the following discrete inf-sup condition:  for any $s\in (1,\infty)$ there exists a constant $\beta_s^{}>0$, independent of $h$,  such that for all $q_{H}^{}\in \calQ_{H}^{}$ the following inequality holds:
\begin{equation}\label{inf-sup-help}
\sup_{\bv_h^{}\in \bV_h^{}\setminus\{\boldsymbol{0}\}}\frac{(q_{H}^{},\dive\bv_h^{})_\Omega^{}}
{|\bv_h^{}|_{1,s',\Omega}^{}}\,\ge\, \beta_s^{}\|q_{H}^{}\|_{0,s,\Omega}^{}.
\end{equation}
The proof of this result, to the best of our knowledge, has not been given previously and thus we report it in the Appendix. It is based on the construction of a Fortin operator $\calI:W^{1,s'}_0(\Omega)^d\to\bV_h^{}$ satisfying
\begin{alignat}{2}
\big(q_H^{},\dive (\bv-\calI(\bv))\big)_\Omega^{} &= 0\qquad &&\textrm{for all}\; q_H^{}\in\calQ_H^{}\quad\textrm{and all}\;\bv\in W^{1,s'}_0(\Omega)^d ,\label{Fortin-1}\,\\
\calI\bv &\to \bv \qquad &&\textrm{strongly in}\;  W^{1,s'}_0(\Omega)^d \;\textrm{as}\; h\to 0. \label{Fortin-2}
\end{alignat}
In addition, \eqref{inf-sup-help} guarantees the existence of a non-trivial subspace of discretely divergence-free functions
\begin{equation}\label{disc-divfree-space}
\bV_{h,\dive}^{}:= \{\bv_h^{}\in\bV_h^{}: (q_H^{},\dive \bv_h^{})_\Omega^{} = 0\;\textrm{ for all}\; q_H^{}\in\calQ_H^{}\}.
\end{equation}


\subsection{Results linked to the discrete Lipschitz truncation}\label{Sec:Lipschitz-Truncation}

In the convergence proof given below we will need the following two results. These are known as
discrete Lipschitz truncation and divergence-free discrete Lipschitz truncation, respectively.
Their proofs are omitted since they are essentially a  rewriting of Corollary~17 and  the proof on pages 1006--1007 in \cite{DKS13} (see also \cite[Lemmas~2.29 and 2.30]{Tabea}).

\begin{lemma}\label{Lemma2.29}
Let $s\in(1,\infty)$. Let us suppose that $\bv_h^{}\in\bV_h^{}$ for all $h>0$ and  $\bv_h^{}\weakcv \boldsymbol{0}$ weakly in $W^{1,s}_0(\Omega)^d$ as $h\to 0$. Then, there exist
\begin{itemize}
\item a double sequence $\{\lambda_{h,j}^{}\}_{h>0,j\in\bN}^{}\subseteq\bR$ such that $\lambda_{h,j}^{}\in [2^{2^j},2^{2^{j+1}-1}]$ for all $h>0,j\in\bN$;
\item a double sequence of open sets $\calB_{h,j}^{}\subseteq\Omega, h>0,j\in\bN$, of the form
\begin{equation}\label{Bhj-set}
\calB_{h,j}^{}= {\rm int}\big(\cup\{ K: K\in\calT_{h,j}^{}\}\big),
\end{equation}
where $\calT_{h,j}^{}$ denotes the collection of some elements of the mesh $\calT_{h}^{}$;
\item a double sequence $\{\bv_{h,j}^{}\}_{h>0,j\in\bN}^{}\subset W^{1,\infty}_0(\Omega)^d$ with $\bv_{h,j}^{}\in\bV_h^{}$ for all
$j\in\bN$ and all $h>0$;
\end{itemize}
satisfying
\begin{itemize}
\item[i.] $\bv_{h,j}^{}=\bv_h^{}$ in $\Omega\setminus\calB_{h,j}^{}$ for all $j\in\bN$ and all $h>0$;
\item[ii.]  there exists a $c(s)>0$ such that
\begin{equation}\label{L2.29ii}
\|\lambda_{h,j}^{}\mathds{1}_{\mathcal{B}_{h,j}^{}}\|_{0,s,\Omega}^{}\le c(s)\,2^{-\frac{j}{s}}\qquad\,\forall\; h>0\;,\; j\in\bN;
\end{equation}
\item[iii.] there exists a $c(s)>0$ such that
\begin{equation}\label{L2.29iii}
\|\nabla\bv_{h,j}^{}\|_{0,\infty,\Omega}^{}\le c(s)\lambda_{h,j}^{}\qquad\,\forall\; h>0\;,\; j\in\bN;
\end{equation}
\item[iv.] for any fixed $j\in\bN$,
\begin{equation}\label{L2.29iv}
\bv_{h,j}^{}\to \boldsymbol{0}\;\textrm{strongly in}\; L^\infty(\Omega)^d\quad\textrm{and}\quad
\nabla \bv_{h,j}^{}\weakcv \boldsymbol{0}\;\textrm{weakly-* in}\; L^\infty(\Omega)^{d\times d},
\end{equation}
as $h\to 0$.
\end{itemize}
\end{lemma}

\begin{lemma}\label{Lemma2.30}
Let $s\in (1,\infty)$ and assume that Assumption~{\rm (A1)} is satisfied. Let $\{\bv_h^{}\}_{h>0}^{}$ be a sequence such that 
$\bv_h^{}\in \bV_{h,\dive}^{}$  for all $h>0$ and such that $\bv_h^{}\weakcv \boldsymbol{0}$ weakly in $W^{1,s}_0(\Omega)^d$ 
as $h\to 0$. Furthermore, let $\{\bv_{h,j}^{}\}_{h>0,j\in\bN}^{}$ be the sequence of Lipschitz truncations given by Lemma~\ref{Lemma2.29}. Then, there
exists a double sequence $\{\bw_{h,j}^{}\}_{h>0,j\in\bN}^{}$ such that
\begin{itemize}
\item[i.] $\bw_{h,j}^{}\in \bV_{h,\dive}^{}$  for all $h>0$ and all $j\in\bN$;
\item[ii.] there exists a $c(s)$ such that
\begin{equation}\label{Lemma2.30ii}
\|\bv_{h,j}^{}-\bw_{h,j}^{}\|_{1,s,\Omega}^{} \le c(s)2^{-\frac{j}{s}}\qquad\forall\, h>0, \; \forall\, j\in\bN;
\end{equation}
\item[iii.] for any fixed $j\in\bN$ the following convergences hold (up to a subsequence, if necessary):
\begin{equation}\label{Lemma2.30iii}
\bw_{h,j}^{}\to\boldsymbol{0}\;\textrm{strongly in}\; L^t(\Omega)^d\quad\textrm{and}\quad \nabla\bw_{h,j}^{}\weakcv 
\boldsymbol{0} \;\textrm{weakly in}\; W^{1,t}_0(\Omega)^{d\times d},
\end{equation}
as $h\to 0$, for all $t<+\infty$.
\end{itemize}
\end{lemma}


\subsection{The stabilising bilinear form and the lifting operator}

The finite element method studied in this work is based on the pair $\bV_h^{}\times\calQ_h^{}$. Since this pair
is not inf-sup stable some form of stabilisation is needed. In this work our proposal is to use the following stabilising bilinear form
\begin{equation}\label{s-definition}
s(q_h^{},t_h^{}) = \sumMTH\sum_{F\in\calF_I^{}(M)}\tau_F^{}\,(\llbracket q_h^{}\rrbracket, \llbracket t_h^{}\rrbracket)_F^{},
\end{equation}
where the stabilisation parameter $\tau_F^{}$ is defined as follows:
\begin{equation}\label{tauF-non-Newtonian}
\tau_F^{}=h_F^{\alpha(r)}\qquad\textrm{where}\qquad
\alpha(r) := \left\{ 
\begin{array}{rll}  1 &&\textrm{if}\; r\ge 2, \\
 1-d+\dfrac{2d}{\tilde{r}} &&\textrm{if}\; r\in \left[\frac{3d}{d+2}, 2\right),\\ ~\\
 1-d+\dfrac{2d}{\tilde{r}_{\rm max}^{}} &= \dfrac{d-1}{3}&\textrm{if}\; r\in \left(\frac{2d}{d+2}, \frac{3d}{d+2}\right).
\end{array}\right.
\end{equation}

The behaviour of $\alpha(r)$ is depicted in Figure~\ref{Fig1}. It can be observed there that the stabilisation gets stronger as $r\to 1$. The reason for this behaviour  will become
clear when we perform the convergence analysis in Section~\ref{Sec:Convergence}.

\begin{remark}\label{Rem:bound-s-average} 
Thanks to $\hypothesis$ and the inf-sup condition \eqref{inf-sup-help} it can be expected to have stability of $\Pi_H^{}(p_h^{})$ (where $p_h^{}$ is the finite element approximation of the pressure). 
The stabilisation is then  built with the aim of controlling $p_h^{}-\Pi_H^{}(p_h^{})$. More precisely, 
using  \eqref{average-bound}, \eqref{bound-lp} (or the inverse inequality \eqref{inverse}),  and the definition of the bilinear form $s(\cdot,\cdot)$ we see that there exists a constant $C>0$ such that
\begin{equation}\label{bound-s}
\|q_h^{}-\Pi_H^{}(q_h^{})\|_{0,s,\Omega}^{}\le C\,h^{\chi}\,s(q_h^{},q_h^{})^{\frac{1}{2}},
\end{equation}
for all $q_h^{}\in\calQ_h^{}$, where 
\begin{equation}
\chi=\left\{ \begin{array}{rl}
\frac{1-\alpha(r)}{2} & \textrm{if}\; s\le 2,\\
\frac{1-d+\frac{2d}{s}-\alpha(r)}{2} & \textrm{if}\; s> 2.
\end{array}\right.
\end{equation}
 It will be  useful in what follows to observe that for $s=\tilde{r}$ we have $\chi\ge 0$.
\end{remark}

Another important ingredient in the definition of the method is a lifting of the pressure jumps defined
with the help of the lowest order Raviart--Thomas basis functions. To define this,
for each  $F\in\calF$ we choose a unique normal vector $\nn_F^{}$. Its orientation is of no importance,
but it needs to point outwards of $\Omega$ if $F\subset \partial\Omega$. Moreover, for each $K\in\calT_h^{}$
such that $F\in\calF_K^{}$, we denote the node in $K$ opposite  $F$ by $\bx_F^{}$.
Using this unique normal vector,
we introduce the lowest order Raviart--Thomas basis function $\bvarphi_F^{}$ defined as
\begin{equation}\label{RT-basis}
\bvarphi_F^{}(\bx)|_K^{}:=\pm \frac{|F|}{d|K|}(\bx-\bx_F^{}),
\end{equation}
and extended by zero outside $\omega_F^{}$. 
In this definition, the sign of the function $\bvarphi_F^{}$ depends on whether the normal vector $\nn_F^{}$
points in or out of $K$. Thanks to its definition,  $\bvarphi_F^{}$ satisfies the following: for every $F'\in\calF$ the
normal component of $\bvarphi_F^{}$ is given by (with the obvious abuse of notation considering that $\nn_{F'}^{}$ is not defined at the boundary of
$F'$):
\begin{equation}\label{dof-RT}
\bvarphi_F^{}\cdot\nn_{F'}^{}=\left\{\begin{array}{ll} 1 & \textrm{if}\; F'=F,\\ 0  & \textrm{otherwise}.\end{array}\right.
\end{equation}

With the help of these Raviart--Thomas basis functions, we define the following operator, which will be
fundamental in the definition of the finite element method:
\begin{gather}
\calL:W^{1,r}(\Omega)^d\times \calQ_h^{}\to W^r(\dive\!; \Omega)\nonumber\\ 
(\bv,q_h^{})\mapsto \calL(\bv,q_h^{}) := \bv + \sumMTH\sum_{F\in\calF_I^{}(M)}\tau_F^{}\llbracket q_h^{}\rrbracket \bvarphi_F^{}.\label{L-definition}
\end{gather}

Since the velocity $\bu$ is bounded in $W^{1,r}_0(\Omega)^d$, then it is bounded in $L^{2\tilde{r}}(\Omega)^d$ as well. In the finite element method
proposed in Section~\ref{Sec:FEM}, we will consider a modified velocity built with the help of the mapping $\calL$ just defined. The following result states that
the stability just mentioned is preserved by the operator $\calL$.

\begin{lemma}\label{Lem:Cont-L} 
There exists a constant $C>0$, independent of $h$, such that
\begin{equation}\label{cont-L}
\|\calL(\bv,q_h^{})\|_{0,2\tilde{r},\Omega}^{}\le C\left\{ |\bv|_{1,r,\Omega}^{}+
s(q_h^{},q_h^{})^{\frac{1}{2}}
\right\},
\end{equation}
for all $(\bv,q_h^{})\in W^{1,r}_0(\Omega)^d\times \calQ_h^{}$.
\end{lemma}

\begin{proof}
Thanks to the embedding \eqref{Sobolev-emb} and denoting
\begin{equation}\label{unc-def}
\bu_{nc}^{}:=\sumMTH\sum_{F\in\calF_I^{}(M)}\tau_F^{}\llbracket q_h^{}\rrbracket \bvarphi_F^{},
\end{equation}
  the following bound follows
\begin{equation}\label{bound-L-1}
\|\calL(\bv,q_h^{})\|_{0,2\tilde{r},\Omega}^{}\le C\,|\bv|_{1,r,\Omega}^{}+\|\bu_{nc}^{}\|_{0,2\tilde{r},\Omega}^{}.
\end{equation}
To bound the second term on the right-hand side of this inequality  we start by noticing that the 
definition of $\bvarphi_F^{}$ (cf. \eqref{RT-basis})  gives
$\|\bvarphi_F^{}\|_{0,\infty,K}^{}\le C$ for each $K$ such that $F\in \calF_K^{}$.
So, let  $K\in \calT_h^{}$ and let $M\in\calT_H^{}$ be the unique macro-element such that $K\subset M$. Then, using the mesh regularity and
the Cauchy-Schwarz inequality we get 
\begin{align*}
\|\bu_{nc}^{}\|_{0,K}^{} 
&=\, \left\| \sum_{F\in\calF_K^{}\cap\calF_I^{}(M)}\tau_F^{}\llbracket p_h^{}\rrbracket \bvarphi_F^{}\right\|_{0,K}^{} \\
&\le\, \sum_{F\in\calF_K^{}\cap\calF_I^{}(M)}\tau_F^{}|\llbracket p_h^{}\rrbracket|\, \|\bvarphi_F^{}\|_{0,K}^{} \\
&\le\, C\,\sum_{F\in\calF_K^{}\cap\calF_I^{}(M)}\tau_F^{} h_F^{\frac{d}{2}}|\llbracket p_h^{}\rrbracket| \\
& \le \, C\,\sum_{F\in\calF_K^{}\cap\calF_I^{}(M)}\tau_F^{} h_F^{1-\frac{d}{2}}(1,|\llbracket p_h^{}\rrbracket|)_F^{} \\
&\le\, C\,h_F^{\frac{1}{2}}\sum_{F\in\calF_K^{}\cap\calF_I^{}(M)}\tau_F^{} \|\llbracket p_h^{}\rrbracket\|_{0,F}^{}.
\end{align*}
Hence, squaring, summing over all the elements, and using the mesh regularity gives
\begin{equation}\label{bound-L-2}
\|\bu_{nc}^{}\|_{0,\Omega}^{} = \left\{\sum_{K\in\calT_h^{}}\|\bu_{nc}^{}\|_{0,K}^2\right\}^{\frac{1}{2}}
\le C\,h^{\frac{1+\alpha(r)}{2}}
\left\{\sumMTH\sum_{F\in\calF_I^{}(M)}\tau_F^{}\|\llbracket p_h^{}\rrbracket \|^2_{0,F}\right\}^{\frac{1}{2}}.
\end{equation}
Thus, using the inverse inequality \eqref{inverse} we arrive at
\begin{equation}
\|\bu_{nc}^{}\|_{0,2\tilde{r},\Omega}^{} \le Ch^{\frac{d(1-\tilde{r})}{2\tilde{r}}}\|\bu_{nc}^{}\|_{0,\Omega}^{}
\le Ch^{\frac{d(1-\tilde{r})+\tilde{r}+\alpha(r)\tilde{r}}{2\tilde{r}}}
\, s(p_h^{},p_h^{})^{\frac{1}{2}}.
\label{bound-L-3}
\end{equation}
To complete the proof we only need to make sure that the exponent of $h$ in \eqref{bound-L-3} is not negative.
Let  $\xi:= d(1-\tilde{r})+\tilde{r}+\alpha(r)\tilde{r}$. If $r\ge 2$ then $\alpha(r)=1$ and $\tilde{r}\le 2$. So,
$\xi=d(1-\tilde{r})+2\tilde{r}= d+(2-d)\tilde{r}\ge 4-d\ge 1$. If $r<2$ then $\alpha(r)\ge \frac{d-1}{3}$ and so
\begin{equation*}
\xi \ge d(1-\tilde{r})+\tilde{r}+\frac{d-1}{3}\tilde{r} = d + \frac{2(1-d)}{3}\tilde{r} \ge d + \frac{2(1-d)}{3}\tilde{r}_{\rm max}^{} = 0.
\end{equation*}
Since in the whole range of values for $r$ we have $\xi\ge 0$, the proof is complete.
\end{proof}


\section{The finite element method}\label{Sec:FEM}

The finite element method studied in this work reads as follows:
find $(\bu_h^{}, p_h^{})\in \bV_h^{}\times \calQ_h^{}$ such that
\begin{align}
(|\nabla\bu_h^{}|^{r-2}\nabla\bu_h^{},\nabla\bv_h^{})_\Omega^{}-(\calL(\bu_h^{},p_h^{})\otimes\bu_h^{},\nabla\bv_h^{})_\Omega^{}-
(p_h^{},\dive\bv_h^{})_\Omega^{}
 &= \langle \ff,\bv_h^{}\rangle_\Omega^{}, \label{discrete-p-laplacian-1}\\
(q_h^{},\dive\bu_h^{})_\Omega^{} 
+ s(p_h^{},q_h^{})&= 0, \label{discrete-p-laplacian-2}
\end{align}
for all $(\bv_h^{},q_h^{})\in \bV_h^{}\times\calQ_h^{}$, 
where $\calL$ is  defined by \eqref{L-definition} and the stabilising bilinear form $s(\cdot,\cdot)$ is defined in \eqref{s-definition}.

\begin{remark} 

\noindent (i) The main differences between \eqref{discrete-p-laplacian-1}, \eqref{discrete-p-laplacian-2} and a standard 
Galerkin method are twofold: first, the stabilising
term involving the jumps of the discrete pressure are added to the formulation to compensate for the fact
that the pair $\bV_h^{}\times \calQ_h^{}$ does not satisfy the discrete inf-sup condition. Additionally,
and perhaps more significantly, the convection velocity $\bu_h^{}$ has been replaced by the modified version $\calL(\bu_h^{},p_h^{})$.
In Lemma~\ref{Lem:L-div-free} this modified velocity will be proved to be solenoidal,  which allows us to analyse the finite element method without the need to rewrite the convection
term in a skew-symmetric form. This will  lead to a convergence result valid in the whole range $r> \frac{2d}{d+2}$.

\noindent (ii) As can be expected, the power of $h$ in the stabilisation parameter depends strongly on the value of $r$. Two important remarks are in order:
\begin{itemize}
\item $\alpha(r)=1$ for all $r\ge 2$;
\item for all $r<2$ we have $\frac{d-1}{3}\le \alpha(r)<1$. 
\end{itemize}
Thus, there is always a positive power of $h$ multiplying the jump terms of the pressure involved in the definition of $s(\cdot,\cdot)$ and
$\calL(\bu_h^{},p_h^{})$, but the stabilisation becomes stronger as $r\to 1$. 
\end{remark}

\subsection{Existence of a solution and a priori bounds}\label{Sec:Existence}
Before exploring the stability  of the scheme, we present the following a priori result concerning
qualitative properties of $\bu_h^{}$ and $\calL(\bu_h^{},p_h^{})$, whenever $(\bu_h^{},p_h^{})$ solves 
\eqref{discrete-p-laplacian-1}, \eqref{discrete-p-laplacian-2}.

\begin{lemma}\label{Lem:L-div-free}
Let $(\bu_h^{},p_h^{})\in\bV_h^{}\times\calQ_h^{}$ be any solution of 
\eqref{discrete-p-laplacian-1}, \eqref{discrete-p-laplacian-2}. Then,

\begin{itemize}
\item[(i)] $\bu_h^{}$ is discretely divergence-free with respect to the coarse space $\calQ_H^{}$, that is,
\begin{equation}\label{div-div-free-H}
(q_H^{},\dive \bu_h^{})_\Omega^{}=0\qquad\forall\, q_H^{}\in \calQ_H^{}.
\end{equation}

\item[(ii)]   $\calL(\bu_h^{},p_h^{})\cdot\nn=0$ on $\partial\Omega$, and
\begin{equation}\label{L-div-free}
\dive\calL(\bu_h^{},p_h^{})=0\quad\textrm{in}\;\Omega.
\end{equation}
\end{itemize}
\end{lemma}

\proof  The proof of (i) is a consequence of the fact that the stabilisation $s(\cdot,\cdot)$ vanishes on the coarse space $\calQ_H^{}$,
that is, $s(q_h^{}, q_H^{})=0$ for all $q_h^{}\in \calQ_h^{}$ and all $q_H^{}\in\calQ_H^{}$.
For (ii), we can follow similar arguments as those presented in \cite[Lemma~3.8]{BV10} and \cite[Lemma~3]{BV11}
(see also \cite[Theorem~3]{ABN18} for a different proof).
$\qed$

The following result states the existence of a solution to the discrete problem  \eqref{discrete-p-laplacian-1}, \eqref{discrete-p-laplacian-2}. In addition, 
it provides uniform a priori  bounds for the sequence of solutions as $h\to 0$.

\begin{theorem} There exists a solution $(\bu_h^{},p_h^{})\in \bV_h^{}\times\calQ_h^{}$  of
 \eqref{discrete-p-laplacian-1}, \eqref{discrete-p-laplacian-2}. 
Moreover,  every solution satisfies the following a priori bound:
\begin{equation}\label{a-priori-bound-p}
|\bu_h^{}|_{1,r,\Omega}^{r}+\|\calL(\bu_h^{},p_h)\|_{0,2\tilde{r},\Omega}^{}+s(p_h^{},p_h^{})
+ \|p_h^{}\|_{0,\tilde{r},\Omega}^{}\le M,
\end{equation}
where $M$ does not depend on $h$.
\end{theorem}

\begin{proof}  The existence of a solution is proved using the argument  used in \cite{GR86} for the Navier--Stokes equation.
First, if $\ff=\boldsymbol{0}$, then $\bu_h^{}=\boldsymbol{0}$ and $p_h^{}=0$ trivially solve \eqref{discrete-p-laplacian-1}, 
\eqref{discrete-p-laplacian-2}. So, we suppose that $\ff\not=\boldsymbol{0}$.
The subspace of $\bV_h^{}\times\calQ_h^{}$ where solutions of \eqref{discrete-p-laplacian-1}, \eqref{discrete-p-laplacian-2} are to 
be sought is given by
\begin{equation}\label{def-Xh}
\bX_h^{}:=\{ (\bv_h^{},q_h^{})\in \bV_h^{}\times\calQ_h^{}\,:\, \dive \calL(\bv_h^{},q_h^{})=0\;\textrm{in}\;\Omega\}.
\end{equation}
Let $T:\bX_h^{}\to [\bX_h^{}]'$ be the mapping defined by
\begin{align}
[T(\bv_h^{},q_h^{}),(\bw_h^{},t_h^{})] &= (|\nabla\bv_h^{}|^{r-2}\nabla\bv_h^{},\nabla\bw_h^{})_\Omega^{} -(\calL(\bv_h^{},q_h^{})\otimes\bv_h^{},\nabla\bw_h^{})_\Omega^{}-
(q_h^{},\dive\bw_h^{})_\Omega^{} \nonumber\\
&\quad +\,(t_h^{},\dive\bv_h^{})_\Omega^{} + s( q_h^{}, t_h^{}) - \langle \ff,\bw_h^{}\rangle_\Omega^{},
\end{align}
that is, the mapping associated with the residual of \eqref{discrete-p-laplacian-1}, \eqref{discrete-p-laplacian-2}. 
For any $(\bv_h^{},q_h^{})\in\bX_h^{}$,  integration by parts gives
\begin{equation}\label{stab-1}
(\calL(\bv_h^{},q_h^{})\otimes\bv_h^{},\nabla\bv_h^{})_\Omega^{} = 0,
\end{equation}
and then, for any $(\bv_h^{},q_h^{})\in\bX_h^{}$,  Young's inequality yields
\begin{align}
[T(\bv_h^{},q_h^{}),(\bv_h^{},q_h^{})] &\ge |\bv_h^{}|^r_{1,r,\Omega}+ s(q_h^{},q_h^{})
-\|\ff\|_{-1,r',\Omega}^{}|\bv_h^{}|_{1,r,\Omega}^{} \nonumber \\
&\ge \frac{1}{r'}|\bv_h^{}|^r_{1,r,\Omega}+ s(q_h^{},q_h^{}) -\frac{1}{r'}\|\ff\|_{-1,r',\Omega}^{r'}. \label{stab-2}
\end{align}
This implies that, for any $(\bv_h^{},q_h^{})\in\bX_h^{}$ such that
\begin{equation}\label{stab-3}
\frac{1}{r'}|\bv_h^{}|^r_{1,r,\Omega}+ s(q_h^{},q_h^{}) = \|\ff\|_{-1,r',\Omega}^{r'},
\end{equation}
we have $[T(\bv_h^{},q_h^{}),(\bv_h^{},q_h^{})]>0$. Thus, using a consequence of Brouwer's fixed point theorem
(see \cite[Ch.~IV, Corollary~1.1]{GR86}) there exists a $(\bu_h^{},p_h^{})\in\bX_h^{}$ 
such that $T(\bu_h^{},p_h^{})=\boldsymbol{0}$. In other words, $(\bu_h^{},p_h^{})\in\bX_h^{}$ solves \eqref{discrete-p-laplacian-1},
\eqref{discrete-p-laplacian-2}. 

In order to prove the a priori bound \eqref{a-priori-bound-p}, we first take $(\bv_h^{},q_h^{})=(\bu_h^{},p_h^{})$  in \eqref{discrete-p-laplacian-1}, \eqref{discrete-p-laplacian-2} and use 
the fact that $\calL(\bu_h^{},p_h^{})$ is solenoidal to arrive at
\begin{equation}\label{a-priori-bound-1}
|\bu_h^{}|_{1,r,\Omega}^{r}+s(p_h^{},p_h^{})\le C\,\|\ff\|_{-1,r',\Omega}^{r'},
\end{equation}
where $C>0$ depends only on $r$. Moreover, the bound on  $\|\calL(\bu_h^{},p_h^{})\|_{0,2\tilde{r},\Omega}^{}$ follows from Lemma~\ref{Lem:Cont-L}
and \eqref{a-priori-bound-1}.

To bound  $\|p_h^{}\|_{0,\tilde{r},\Omega}^{}$ we consider the projection $\Pi_H^{}$ defined in Section~\ref{FEspaces} and write
\begin{equation}\label{first-p-bound}
\|p_h^{}\|_{0,\tilde{r},\Omega}^{}\le \|p_h^{}-\Pi_H^{}(p_h^{})\|_{0,\tilde{r},\Omega}^{} + \|\Pi_H^{}(p_h^{})\|_{0,\tilde{r},\Omega}^{}.
\end{equation}
First, using the result stated in Remark~\ref{Rem:bound-s-average} 
(that is, \eqref{bound-s} with $s=\tilde{r}$) to deduce that
\begin{equation}\label{p-tres}
\|\Pi_{H}^{}(p_h^{})-p_h^{}\|_{0,\tilde{r},\Omega}^{}
\le C\,h^{\chi}s(p_h^{},p_h^{})^{\frac{1}{2}}\le C,
\end{equation}
where $\chi\ge 0$. Next, since $\Pi_H^{}(p_h^{})\in\calQ_H^{}$, thanks to \eqref{inf-sup-help}  there exists a $\tilde{\bw}_h^{}\in\bV_h^{}$ such that
$|\tilde{\bw}_h^{}|_{1,\tilde{r}',\Omega}^{}=1$ and 
\begin{align}
&\beta_r^{}\|\Pi_{H}^{}(p_h^{})\|_{0,\tilde{r},\Omega}^{} \,\le\, (\Pi_{H}^{}(p_h^{}), \dive\tilde{\bw}_h^{})_\Omega^{}
\nonumber\\
&=\, (\Pi_{H}^{}(p_h^{})-p_h^{}, \dive\tilde{\bw}_h^{})_\Omega^{} + 
 (p_h^{}, \dive\tilde{\bw}_h^{})_\Omega^{}\nonumber\\
&=\, (\Pi_{H}^{}(p_h^{})-p_h^{}, \dive\tilde{\bw}_h^{})_\Omega^{} + 
(|\nabla\bu_h^{}|^{r-2}\nabla\bu_h^{},\nabla\tilde{\bw}_h^{})_\Omega^{} -
(\calL(\bu_h^{},p_h^{})\otimes \bu_h^{},\nabla\tilde{\bw}_h^{})_\Omega^{}-\langle \ff,\tilde{\bw}_h^{}\rangle_\Omega^{}\nonumber\\
&= I+ II + III + IV,
\label{p-uno}
\end{align}
where we have also used  that $(\bu_h^{},p_h^{})$ solves \eqref{discrete-p-laplacian-1}, \eqref{discrete-p-laplacian-2}.
The bounds for the above terms proceed  using  H\"older's inequality, 
 $\tilde{r}\le r'$ (and then $r\le \tilde{r}'$), $|\tilde{\bw}_h^{}|_{1,\tilde{r}',\Omega}^{}=1$,  \eqref{p-tres}, \eqref{Sobolev-emb} , the bound for $\|\calL(\bu_h^{},p_h^{})\|_{0,2\tilde{r},\Omega}^{}$, and 
\eqref{a-priori-bound-1} as follows:
\begin{align}
I &\le \|\Pi_{H}^{}(p_h^{})-p_h^{}\|_{0,\tilde{r},\Omega}^{}|\tilde{\bw}_h^{}|_{1,\tilde{r}',\Omega}^{}
\le C, \label{p-tres-bis}\\
II&\le\, \left(\int_\Omega|\nabla\bu_h^{}|^{(r-1)\tilde{r}}\right)^{\frac{1}{\tilde{r}}}|\tilde{\bw}_h^{}|_{1,\tilde{r}',\Omega}^{}
\le C\, \left(\int_\Omega|\nabla\bu_h^{}|^{(r-1)r'}\right)^{\frac{1}{r'}} = C\, |\bu_h^{}|_{1,r,\Omega}^{r-1}\le C, \label{p-cinco}\\
III&\le \|\calL(\bu_h^{},p_h^{})\|_{0,2\tilde{r},\Omega}^{}\|\bu_h^{}\|_{0,2\tilde{r},\Omega}^{}
|\tilde{\bw}_h^{}|_{1,\tilde{r}',\Omega}^{} \le 
C\, \|\calL(\bu_h^{},p_h^{})\|_{0,2\tilde{r},\Omega}^{}|\bu_h^{}|_{1,r,\Omega}^{} \le C, \label{p-seis} \\
IV &\le \|\ff\|_{-1,r',\Omega}^{}|\tilde{\bw}_h^{}|_{1,r,\Omega}^{}\le C\,\|\ff\|_{-1,r',\Omega}^{}\,|\tilde{\bw}_h^{}|_{1,\tilde{r}',\Omega}^{}=C\,\|\ff\|_{-1,r',\Omega}^{}.
\label{p-siete}
\end{align}
Thus, the proof follows by inserting the above bounds on $I,\ldots,IV$ and \eqref{p-tres}  in \eqref{first-p-bound}.
\end{proof}

%

\section{Convergence to a weak solution}\label{Sec:Convergence}

In this section we analyse the convergence of the finite element scheme \eqref{discrete-p-laplacian-1}, \eqref{discrete-p-laplacian-2}. The convergence proof is divided into two
cases in order to distinguish between the situations when a solution $\bu$ of \eqref{r-lap-1}, \eqref{r-lap-2} can, and cannot, be used as a test function in \eqref{r-lap-1}. 

\begin{theorem} Let $r\in (\frac{2d}{d+2},\infty)$. Then, there exists a subsequence, still denoted by $(\bu_h^{},p_h^{})$, such that
\begin{align}
&\bu_h^{} \rightharpoonup \bu\quad\textrm{weakly in}\; W^{1,r}_0(\Omega)^d;  \label{weak_u-r}\\
&\bu_h^{} \to \bu\quad\textrm{strongly in}\; L^s(\Omega)^d\quad\textrm{for}\; s\in\, [1,2\tilde{r}); 
\label{strong_u-r}\\
&\calL(\bu_h^{},p_h^{}) \to \bu\quad\textrm{strongly in}\; L^s(\Omega)^d,\;
\textrm{for}\;s\in\, [1,2\tilde{r}); 
 \label{strong_L(u,p)-r}\\
&p_h^{} \rightharpoonup p \quad\textrm{weakly in}\; L^{\tilde{r}}(\Omega); \label{weak_p-r} \\
& 
\textrm{if}\; r\ge \frac{3d}{d+2},\quad\textrm{then}\;
s(p_h^{},p_h^{})\to 0\qquad\textrm{as}\; h\to 0. \label{strong-jump-r}
\end{align}
In addition, $p\in L_0^{\tilde{r}}(\Omega)$, and $(\bu,p)$  solves \eqref{r-lap-1}, \eqref{r-lap-2}.
\end{theorem}

\begin{proof}
The proofs of \eqref{weak_u-r} and \eqref{weak_p-r} follow using \eqref{a-priori-bound-p} and the reflexivity of 
$W^{1,r}_0(\Omega)$ and $L^{\tilde{r}}(\Omega)$ for $r\in (1,\infty)$. In addition, $p$ has zero average since
\begin{equation}\label{p-zero-average}
(p,1)_\Omega^{}=\lim_{h\to 0}(p_h^{},1)_\Omega^{}=0.
\end{equation}
The proof of \eqref{strong_u-r} is a consequence of the Rellich--Kondrachov Theorem (see, e.g.
\cite[Theorem~9.16]{Brezis}).  Moreover, the bound \eqref{bound-L-3} implies that for every $s< 2\tilde{r}$ there exists a number $\xi>0$ such that
\begin{equation}\label{unc-strong-to-zero}
\|\bu_{nc}^{}\|_{0,s,\Omega}^{}\le Ch^{\xi} s(p_h^{},p_h^{})^{\frac{1}{2}},
\end{equation}
 so \eqref{a-priori-bound-p} yields $\|\bu_{nc}^{}\|_{0,s,\Omega}^{}\to 0$ as $ h\to 0$ for all $s< 2\tilde{r}$. Together with   \eqref{strong_u-r} this proves
\eqref{strong_L(u,p)-r}.


We now start the process of identifying the partial differential equation satisfied by the limits $\bu$ and $p$.
Let $\bv\in C_0^{\infty}(\Omega)^d$ be arbitrary, and let $\bv_h^{}\in\bV_h^{}$
be its Scott--Zhang interpolant.  Using \eqref{a-priori-bound-p} we first get that
\begin{equation}\label{extra-1}
\|\,|\nabla\bu_h^{}|^{r-2}\nabla\bu_h^{}\|_{0,r',\Omega}^{r'} = \int_\Omega |\nabla\bu_h^{}|^{\frac{(r-1)r}{(r-1)}} =
\|\nabla\bu_h^{}\|^r_{0,r,\Omega}\le C,
\end{equation}
and thus there exists a $\boldsymbol{S}\in L^{r'}(\Omega)^{d\times d}$ such that (up to a subsequence)
\begin{equation}\label{extra-2}
|\nabla\bu_h^{}|^{r-2}\nabla\bu_h^{}\rightharpoonup \boldsymbol{S}\quad\textrm{weakly in}\; L^{r'}(\Omega)^{d\times d}.
\end{equation}
So, since $\bv_h^{}$ converges to $\bv$ strongly in $W^{1,r}_0(\Omega)$ we have
\begin{equation}\label{uno}
(|\nabla\bu_h^{}|^{r-2}\nabla\bu_h^{},\nabla\bv_h^{})_\Omega^{}\to (\boldsymbol{S},\nabla\bv)_\Omega^{}\quad\textrm{as}\; h\to 0.
\end{equation}

Next, thanks to \eqref{weak_p-r}  and the strong convergence of $\bv_h^{}$ to $\bv$ 
in $W^{1,\tilde{r}'}_0(\Omega)$  (see \eqref{conver-Pih}) the following holds:
\begin{equation}\label{cinco}
(p_h^{},\dive\bv_h^{})_\Omega^{} = (p_h^{},\dive(\bv_h^{}-\bv))_\Omega^{}
+ (p_h^{},\dive \bv)_\Omega^{}\to 0+(p,\dive\bv)_\Omega^{}\qquad\textrm{as}\;h\to 0.
\end{equation}

To treat the convection term, \eqref{strong_u-r} and \eqref{strong_L(u,p)-r} imply that
\begin{equation}\label{conver_of_Lxu}
\calL(\bu_h^{},p_h^{})\otimes \bu_h^{} \to \bu\otimes\bu \qquad\textrm{strongly in}\;L^s(\Omega)^d\quad\textrm{for all}\; s< \tilde{r},
\end{equation}
which, together with the fact that $\bv_h^{}\to \bv$ strongly in $W^{1,s'}_0(\Omega)^d$, prove that
\begin{equation}
(\calL(\bu_h^{},p_h^{})\otimes\bu_h^{},\nabla\bv_h^{})_\Omega^{} \to (\bu\otimes\bu,\nabla\bv)_\Omega^{}\qquad \textrm{as}\; h\to 0. \label{ocho} 
\end{equation}

Thus, $(\boldsymbol{S},\bu,p)$ solves a problem related to \eqref{r-lap-1}. 
In fact, since $(\bu_h^{},p_h^{})$
satisfies \eqref{discrete-p-laplacian-1}, then applying \eqref{uno}, \eqref{cinco}, and \eqref{ocho}, we arrive at
\begin{equation}\label{once}
\begin{array}{ccccccc}
(|\nabla\bu_h^{}|^{r-2}\nabla\bu_h^{},\nabla\bv_h^{})_\Omega^{} & - & (\calL(\bu_h^{},p_h^{})\otimes\bu_h^{},\nabla\bv_h^{})_\Omega^{}
&-& (p_h^{},\dive\bv_h^{})_\Omega^{} &=& \langle\ff,\bv_h^{}\rangle_\Omega^{}  \\
\downarrow & & \downarrow & & \downarrow & & \downarrow \\
(\boldsymbol{S},\nabla\bv)_\Omega^{}  
& - & (\bu\otimes\bu,\nabla\bv)_\Omega^{}
&-& (p,\dive\bv)_\Omega^{} &=& \langle\ff,\bv\rangle_\Omega^{} , 
\end{array}
\end{equation}
as $h\to 0$,  and using the density of $C_0^\infty(\Omega)^d$ in $W^{1,\tilde{r}'}_0(\Omega)^d$, 
$(\bu,p,\boldsymbol{S})$ satisfies
\begin{equation}\label{extra-3}
(\boldsymbol{S},\nabla\bv)_\Omega^{}  - (\bu\otimes\bu,\nabla\bv)_\Omega^{}
- (p,\dive\bv)_\Omega^{} = \langle\ff,\bv\rangle_\Omega^{} \qquad\forall\, \bv\in W^{1,\tilde{r}'}_0(\Omega).
\end{equation}

To show that $\bu$ is solenoidal we consider $q\in C^\infty_0(\Omega)$, integrate by parts, and
use \eqref{strong_L(u,p)-r} and the fact $\calL(\bu_h^{},p_h^{})$ is solenoidal to obtain
\begin{equation}\label{doce}
(\dive\bu,q)_\Omega^{}=-(\bu,\nabla q)_\Omega^{} = -\lim_{h\to 0}(\calL(\bu_h^{},p_h^{}),\nabla q)_\Omega^{}
=\lim_{h\to 0}(\dive\calL(\bu_h^{},p_h^{}),q)_\Omega^{}=0,
\end{equation}
and then $\dive\bu=0$ in the distributional sense.


To prove that that $(\bu,p)$ solves \eqref{r-lap-1}, \eqref{r-lap-2}
it only remains to show that $\boldsymbol{S}=|\nabla\bu|^{r-2}\nabla \bu$. 
The proof  of this  will be split into two cases, labelled (i) and (ii) below.

\medskip

\noindent (i) \underline{$r\ge \frac{3d}{d+2}$ :} 
In this case we use a classical result  commonly referred to as the {\sl Minty trick} (see, e.g., \cite[Lemma~2.13]{Roubicek}).
Let $\bv\in W^{1,r}_0(\Omega)$, and let $\bv_h^{}$ be its Scott--Zhang interpolant. Since the $r$-Laplacian
operator is monotone (see, e.g, \cite{Brezis}) we have
\begin{align}\label{extra-4}
0 &\le (|\nabla\bv|^{r-2}\nabla\bv- |\nabla\bu_h^{}|^{r-2}\nabla\bu_h^{},\nabla(\bv-\bu_h^{}))_\Omega^{} \nonumber\\
&= (|\nabla\bv|^{r-2}\nabla\bv- |\nabla\bu_h^{}|^{r-2}\nabla\bu_h^{},\nabla(\bv-\bv_h^{}))_\Omega^{} \nonumber\\
&\qquad + (|\nabla\bv|^{r-2}\nabla\bv,\nabla(\bv_h^{}-\bu_h^{}))_\Omega^{} - ( |\nabla\bu_h^{}|^{r-2}\nabla\bu_h^{},\nabla(\bv_h^{}-\bu_h^{}))_\Omega^{}\nonumber \\
&= \mathcal{A}+\mathcal{B}+\mathcal{C}.
\end{align}
Using that $\bv_h^{}$ converges strongly to $\bv$ in $W^{1,r}_0(\Omega)^d$ and $\bu_h^{}$ converges weakly to $\bu$ in $W^{1,r}_0(\Omega)^d$, and \eqref{a-priori-bound-p} we easily get
\begin{align}
\begin{aligned}
\mathcal{A}&\le \||\nabla\bv|^{r-2}\nabla\bv- |\nabla\bu_h^{}|^{r-2}\nabla\bu_h^{}\|_{0,r',\Omega}^{}\,|\bv-\bv_h^{}|_{1,r,\Omega}^{}\to 0
\qquad\textrm{as}\; h\to 0, \label{extra-5}\\
\mathcal{B}&= (|\nabla\bv|^{r-2}\nabla\bv,\nabla(\bv_h^{}-\bu_h^{}))_\Omega^{} \to (|\nabla\bv|^{r-2}\nabla\bv,\nabla(\bv-\bu^{}))_\Omega^{}
\qquad\textrm{as}\; h\to 0. 
\end{aligned}
\end{align}
To treat $\mathcal{C}$ we use that $(\bu_h^{},p_h^{})$ solves the discrete problem \eqref{discrete-p-laplacian-1}, \eqref{discrete-p-laplacian-2}, as follows:
\begin{align}
\mathcal{C}&= - ( |\nabla\bu_h^{}|^{r-2}\nabla\bu_h^{},\nabla(\bv_h^{}-\bu_h^{}))_\Omega^{} \nonumber\\
&=  -( |\nabla\bu_h^{}|^{r-2}\nabla\bu_h^{},\nabla\bv_h^{})_\Omega^{} +  
( |\nabla\bu_h^{}|^{r-2}\nabla\bu_h^{},\nabla\bu_h^{})_\Omega^{}  \nonumber\\
&= -( |\nabla\bu_h^{}|^{r-2}\nabla\bu_h^{},\nabla\bv_h^{})_\Omega^{} + 
\underbrace{(\calL(\bu_h^{},p_h^{})\otimes\bu_h^{},\nabla \bu_h^{})_\Omega^{}}_{=0} + (p_h^{},\dive\bu_h^{})_\Omega^{}
+ \langle\ff,\bu_h^{}\rangle_\Omega^{} \nonumber\\
&=  -( |\nabla\bu_h^{}|^{r-2}\nabla\bu_h^{},\nabla\bv_h^{})_\Omega^{} 
- s(p_h^{},p_h^{}) 
+ \langle\ff,\bu_h^{}\rangle_\Omega^{} \nonumber\\
&= 
\mathcal{D}-\mathcal{E}+\mathcal{F}
. \label{extra-6}
\end{align}
Thus, from \eqref{extra-4} we get
\begin{equation*}
0\le 
\mathcal{E} = \mathcal{A} +\mathcal{B}
+ \mathcal{D}+ \mathcal{F},
\end{equation*}
and, taking the limit when $h\to 0$ on both sides of this inequality using that $\bv_h^{}\to\bv$
strongly in $W^{1,r}_0(\Omega)^d$,  $\bu_h^{}\rightharpoonup\bu$ weakly in $W^{1,r}_0(\Omega)^d$,
and $|\nabla\bu_h^{}|^{r-2}\nabla\bu_h^{}\rightharpoonup \boldsymbol{S}$ weakly in $L^{r'}(\Omega)^{d\times d}$, 
we obtain
\begin{equation}\label{extra-7}
0\le \lim_{h\to 0} s(p_h^{},p_h^{}) 
\le (|\nabla\bv|^{r-2}\nabla\bv , \nabla(\bv-\bu))_\Omega^{}-(\boldsymbol{S},\nabla\bv)_\Omega^{}+\langle\ff,\bu\rangle_\Omega^{}.
\end{equation}
It only remains to show that $\langle\ff,\bu\rangle_\Omega^{}= (\boldsymbol{S},\nabla\bu)_\Omega^{}$ to prove
that $\boldsymbol{S}$ and $\bu$ satisfy
\begin{equation}\label{extra-8}
0\le \lim_{h\to 0} s(p_h^{},p_h^{}) 
\le(|\nabla\bv|^{r-2}\nabla\bv-\boldsymbol{S},\nabla(\bv-\bu))_\Omega\qquad\forall\bv\in W^{1,r}_0(\Omega),
\end{equation}
and then the monotonicity of the $r$-Laplacian and
 an application of \cite[Lemma~2.13]{Roubicek} gives $\boldsymbol{S}=|\nabla\bu|^{r-2}\nabla \bu$.
Since $r\ge \frac{3d}{d+2}$ then $r=\tilde{r}'$ and so $\bu\in W^{1,\tilde{r}'}_0(\Omega)$. Hence, by taking 
  $\bv=\bu$ as test function in \eqref{extra-3} we obtain
\begin{equation}\label{extra-10}
\langle\ff,\bu\rangle_\Omega^{} = (\boldsymbol{S},\nabla\bu)_\Omega^{} - 
\underbrace{(\bu\otimes\bu,\nabla\bu)_\Omega^{}}_{=0}- \underbrace{(p,\dive\bu)_\Omega^{}}_{=0}
= (\boldsymbol{S},\nabla\bu)_\Omega^{},
\end{equation}
thus proving that $(\boldsymbol{S},\nabla\bu)_\Omega^{}=\langle\ff,\bu\rangle_\Omega^{}$. Hence $(\bu,p)$ solves the continuous problem 
\eqref{r-lap-1}, \eqref{r-lap-2}. 

Finally, \eqref{strong-jump-r} follows by taking $\bv=\bu$ in \eqref{extra-8}. 


\medskip

\noindent (ii) \underline{ $r \in\left( \frac{2d}{d+2},\frac{3d}{d+2}\right)$ :} 
For this case we are not able to use the fundamental step of taking $\bv=\bu$ as test function in \eqref{extra-3} to conclude $\boldsymbol{S}=|\nabla\bu|^{r-2}\nabla\bu$.
So, we need to appeal to the results concerning discrete Lipschitz truncation described in Section~\ref{Sec:Lipschitz-Truncation}, and use the Minty trick once again.
To conclude that $\boldsymbol{S}=|\nabla\bu|^{r-2}\nabla\bu$ in $\Omega$ we need to show that
\begin{equation}\label{necessary-1}
\lim_{h\to 0}\big(|\nabla\bu_h^{}|^{r-2}\nabla\bu_h^{}-|\nabla\bu|^{r-2}\nabla\bu,\nabla(\bu_h^{}-\bu)\big)_\Omega^{} = 0.
\end{equation}
In fact, using \eqref{necessary-1}  it is simple to prove that
\begin{equation}\label{necessary-2}
\lim_{h\to 0} \big(|\nabla\bu_h^{}|^{r-2}\nabla\bu_h^{},\nabla\bu_h^{})_\Omega^{} = (\boldsymbol{S},\nabla \bu)_\Omega^{},
\end{equation}
which implies that, for all $\bv\in W^{1,r}_0(\Omega)^d$ we have
\begin{align*}
0 &\le \lim_{h\to 0} \big(|\nabla\bu_h^{}|^{r-2}\nabla\bu_h^{}-|\nabla\bv|^{r-2}\nabla\bv,\nabla(\bu_h^{}-\bv)\big)_\Omega^{} \\
&= \big(\boldsymbol{S}-|\nabla\bv|^{r-2}\nabla\bv,\nabla(\bu-\bv)\big)_\Omega^{},
\end{align*}
and the application of the Minty trick gives $\boldsymbol{S}=|\nabla\bu|^{r-2}\nabla\bu$.

To prove \eqref{necessary-1}
let $\textrm{H}_h^{}:= ( |\nabla\bu_h^{}|^{r-2}\nabla\bu_h^{}-|\nabla\bu|^{r-2}\nabla\bu):\nabla(\bu_h^{}-\bu)$. Since the $r$-Laplacian is  monotone,  then $\textrm{H}_h^{}\ge 0$
almost everywhere in $\Omega$ leading to
\begin{equation}\label{Hh-uno}
\liminf_{h\to 0} \int_\Omega \textrm{H}_h^{}(\bx)\textrm{d}\bx \ge 0.
\end{equation}

To prove the converse to \eqref{Hh-uno}, 
let $\bv_h^{}:= \bu_h^{}-\calI(\bu)$, where $\calI$ is the Fortin operator satisfying \eqref{Fortin-1}, \eqref{Fortin-2}. First, $\bu_h^{},\calI(\bu)\in \bV_{h,\dive}^{}$ and $\bv_h^{}\weakcv \boldsymbol{0}$
weakly in $W_0^{1,r}(\Omega)^d$ as $h\to 0$. Let now $\{\bv_{h,j}^{}\}_{h>0,j\in\bN}^{}$ and $\{\bw_{h,j}^{}\}_{h>0,j\in\bN}^{}$ be the sequences defined in Lemmas \ref{Lemma2.29} and \ref{Lemma2.30}, respectively, and let
$\{\calB_{h,j}^{}\}_{h>0,j\in\bN}^{}$ be the sets defined in Lemma~\ref{Lemma2.29}. 
First, thanks to \eqref{a-priori-bound-p}, $\textrm{H}_h^{}$ is uniformly bounded in $L^1(\Omega)$, and then using H\"older's inequality and \eqref{L2.29ii} we get
\begin{align}
\int_\Omega \textrm{H}_h^{\frac{1}{2}}\textrm{d}\bx &= \int_{\calB_{h,j}^{}}\textrm{H}_h^{\frac{1}{2}} \textrm{d}\bx + \int_{\Omega\setminus \calB_{h,j}^{}}\textrm{H}_h^{\frac{1}{2}} \textrm{d}\bx\nonumber\\
&\le |\calB_{h,j}^{}|^{\frac{1}{2}}\left\{ \int_{\calB_{h,j}^{}}\textrm{H}_h^{}\textrm{d}\bx\right\}^{\frac{1}{2}} +   
|\Omega\setminus \calB_{h,j}^{}|^{\frac{1}{2}}\left\{ \int_{\Omega\setminus \calB_{h,j}^{}}\textrm{H}_h^{}\textrm{d}\bx\right\}^{\frac{1}{2}}\nonumber\\
&\le C2^{-\frac{j}{r}}+|\Omega|^{\frac{1}{2}}\,\mathfrak{A}^{\frac{1}{2}}. 
\end{align}
The goal will be to show that $\mathfrak{A}$ is bounded by $C2^{-\frac{j}{r}}$ plus a term that tends to zero with $h$,  ultimately proving that
\begin{equation}
\limsup_{h\to 0}\int_\Omega \textrm{H}_h^{\frac{1}{2}} \textrm{d}\bx\le C2^{-\frac{j}{r}},
\end{equation}
for every $j\in\mathbb{N}$, which combined with \eqref{Hh-uno} will prove \eqref{necessary-1}. 

To bound $\mathfrak{A}$ we start by decomposing the error $\bu_h^{}-\bu$ as $\bu_h^{}-\bu = \bv_h^{}+\calI(\bu)-\bu$, define 
$\calG_h^{}:=|\nabla\bu_h^{}|^{r-2}\nabla\bu_h^{}-|\nabla\bu|^{r-2}\nabla\bu$,
and thus write
\begin{equation}
\mathfrak{A}= \int_{\Omega\setminus \calB_{h,j}^{}}\calG_h^{}:\nabla\bv_h^{}\textrm{d}\bx+\int_{\Omega\setminus \calB_{h,j}^{}}\calG_h^{}:\nabla(\calI(\bu)-\bu) = \mathfrak{B}+\mathfrak{C}.
\end{equation}
Since $\calI(\bu)\to \bu$ strongly in $W^{1,r}_0(\Omega)^d$ and $\calG_h^{}$ is uniformly bounded in $L^{r'}(\Omega)^{d\times d}$
(thanks to \eqref{a-priori-bound-p}), $\mathfrak{C}\to 0$ as $h\to 0$. Moreover, since $\bv_{h,j}^{}=\bv_h^{}$ in $\Omega\setminus\calB_{h,j}^{}$, then
\begin{align}
\mathfrak{B} &= \int_{\Omega\setminus \calB_{h,j}^{}}\calG_h^{}:\nabla\bv_{h,j}^{} \nonumber\\
&= \int_{\Omega}\calG_h^{}:\nabla\bv_{h,j}^{} - \int_{\Omega}\calG_h^{}:\nabla\bv_{h,j}^{} \mathds{1}_{\mathcal{B}_{h,j}^{}} \nonumber\\
&= \int_{\Omega}\calG_h^{}:\nabla(\bv_{h,j}^{}-{\bw}_{h,j}^{})+  \int_{\Omega}\calG_h^{}:\nabla {\bw}_{h,j}^{} - \int_{\Omega}\calG_h^{}:\nabla\bv_{h,j}^{} \mathds{1}_{\mathcal{B}_{h,j}^{}} \nonumber\\
&= \mathfrak{D}+\mathfrak{E}+\mathfrak{F}.
\end{align}
 H\"older's inequality, \eqref{Lemma2.30ii},  \eqref{L2.29iii}, and \eqref{L2.29ii} yield the bounds
\begin{align}
|\mathfrak{D}|&\le \|\calG_h^{}\|_{0,r',\Omega}^{}\|\nabla (\bv_{h,j}^{}- {\bw}_{h,j}^{})\|_{0,r,\Omega}^{}\le C2^{-\frac{j}{r}}, \nonumber\\
|\mathfrak{F}|&\le \|\calG_h^{}\|_{0,r',\Omega}^{}\| \nabla\bv_{h,j}^{} \mathds{1}_{\mathcal{B}_{h,j}^{}}\|_{0,r,\Omega}^{} \le C2^{-\frac{j}{r}},
\end{align}
for all $h>0$. Moreover, $\mathfrak{E}$ is decomposed as follows
\begin{equation}
\mathfrak{E}= \int_\Omega |\nabla\bu_h^{}|^{r-2}\nabla\bu_h^{}:\nabla\bw_{h,j}^{}- \int_\Omega |\nabla\bu|^{r-2}\nabla\bu:\nabla\bw_{h,j}^{} = \mathfrak{G}+\mathfrak{H}.
\end{equation}
Since $\bw_{h,j}^{}\weakcv \boldsymbol{0}$ weakly in $W^{1,r}_0(\Omega)^d$ then $\mathfrak{H}\to 0$ as $h\to 0$. 
The only remaining term to deal with is $\mathfrak{G}$. We start by using that $(\bu_h^{},p_h^{})$ solves
\eqref{discrete-p-laplacian-1}, \eqref{discrete-p-laplacian-2}  to rewrite $\mathfrak{G}$ as follows
\begin{equation}
\mathfrak{G}= (\calL(\bu_h^{},p_h^{})\otimes \bu_h^{}, \nabla\bw_{h,j}^{})_\Omega^{}+ (p_h^{},\dive\bw_{h,j}^{})_\Omega^{}-\langle\boldsymbol{f}, \bw_{h,j}^{}\rangle_\Omega^{} .
\end{equation}
The convection term above is treated as follows: using that for any fixed $s<+\infty$, $\nabla\bw_{h,j}^{}$ is uniformly bounded in $L^s(\Omega)^{d\times d}$, then
 \eqref{conver_of_Lxu} applied to $\hat{s}=\frac{1+\tilde{r}}{2}<\tilde{r}$ yields the bound
\begin{equation}
(\calL(\bu_h^{},p_h^{})\otimes\bu_h^{}-\bu\otimes\bu, \nabla\bw_{h,j}^{})_\Omega^{} \le \| \calL(\bu_h^{},p_h^{})\otimes\bu_h^{}-\bu\otimes\bu\|_{0,\hat{s},\Omega}^{}\|\nabla\bw_{h,j}^{}\|_{0,\hat{s}',\Omega}^{}\to 0,
\end{equation}
as $h\to 0$, and then
\begin{equation}
(\calL(\bu_h^{},p_h^{})\otimes\bu_h^{}, \nabla\bw_{h,j}^{})_\Omega^{}  = (\calL(\bu_h^{},p_h^{})\otimes\bu_h^{}-\bu\otimes\bu, \nabla\bw_{h,j}^{})_\Omega^{}+ 
 (\bu\otimes\bu, \nabla\bw_{h,j}^{})_\Omega^{}\to 0,
\end{equation}
as $h\to 0$. Moreover, $\langle\boldsymbol{f}, \bw_{h,j}^{}\rangle_\Omega^{} \to 0$ as $h\to 0$. Finally, for the remaining term in $\mathfrak{G}$ we get,
 by applying that $\bw_{h,j}^{}\in \bV_{h,\dive}^{}$, Cauchy--Schwarz's inequality and \eqref{bound-s}:
\begin{align}
(p_h^{},\dive\bw_{h,j}^{})_\Omega^{} &= (p_h^{}-\Pi_H^{}(p_h^{}),\dive\bw_{h,j}^{})_\Omega^{} \nonumber\\
&\le \|p_h^{}-\Pi_H^{}(p_h^{})\|_{0,\Omega}^{}\|\dive\bw_{h,j}^{}\|_{0,\Omega}^{} \nonumber\\
&\le Ch^{\frac{1-\alpha(r)}{2}}s(p_h^{},p_h^{})^{\frac{1}{2}}\|\dive\bw_{h,j}^{}\|_{0,\Omega}^{}  \to 0,
\end{align}
as $h\to 0$, since $\alpha(r)=\frac{d-1}{3}<1$ for all $r< \frac{3d}{d+2}$, and $s(p_h^{},p_h^{})$ and $\|\dive\bw_{h,j}^{}\|_{0,\Omega}^{}$ are uniformly bounded in $h$ and $j$. 

Collecting all the above bounds the following can be concluded 
\begin{equation}
\mathfrak{A} = \mathfrak{B}+\mathfrak{C} = \mathfrak{D}+\mathfrak{E}+\mathfrak{F}+\mathfrak{C} \le C2^{-\frac{j}{r}}+ \mathfrak{G}+\mathfrak{H}+\mathfrak{C},
\end{equation}
and since $\mathfrak{G}+\mathfrak{H}+\mathfrak{C}\to 0$ as $h\to 0$ for every fixed $j\in \mathbb{N}$, then, for every $j\in \mathbb{N}$ we get
$\limsup_{h\to 0} \int_\Omega \textrm{H}_h^{\frac{1}{2}}(\bx)\textrm{d}\bx \le C2^{-\frac{j}{r}}$ for every $j\in \mathbb{N}$, and thus
\begin{equation}
\limsup_{h\to 0} \int_\Omega \textrm{H}_h^{\frac{1}{2}}(\bx)\textrm{d}\bx \le 0.
\end{equation}
So, $\int_\Omega \textrm{H}_h^{\frac{1}{2}}\textrm{d}\bx \to 0$, which means that, up to a subsequence if necessary, $\textrm{H}_h^{\frac{1}{2}} \to 0$ almost everywhere in $\Omega$, and thus
$\textrm{H}_h^{} \to 0$ almost everywhere in $\Omega$. This, together with \eqref{Hh-uno}, proves \eqref{necessary-1} and thus $\boldsymbol{S}=|\nabla\bu|^{r-2}\nabla\bu$ almost everywhere in $\Omega$. 
Hence, $(\bu,p)$ solves the continuous problem \eqref{r-lap-1}, \eqref{r-lap-2}. 
\end{proof}


\subsection{Strong convergence}

The convergence results proved in the last section can be strengthened. In fact, in this section we prove that the velocity and pressure converge
strongly, at least for an appropriate range of values of $r$ in the case of the pressure. We start with the proof of the strong convergence of the velocity.

\begin{theorem}\label{Theo:strong-conv-u} For every $r>\frac{2d}{d+2}$ the discrete velocity $\bu_h^{}$ converges to $\bu$ strongly in $W^{1,r}_0(\Omega)^d$.
\end{theorem}

\begin{proof} 
We start by considering the case when $r \ge \frac{3d}{d+2}$.
Using the discrete problem \eqref{discrete-p-laplacian-1},  \eqref{discrete-p-laplacian-2}, \eqref{strong-jump-r},  and \eqref{r-lap-1} with $\bv=\bu$ we get
\begin{align}
\lim_{h\to 0}\big (|\nabla\bu_h^{}|^{r-2}\nabla\bu_h^{},\nabla\bu_h^{}\big)_\Omega^{} &= \lim_{h\to 0}\big\{ \underbrace{(\calL(\bu_h^{},p_h^{})\otimes \bu_h^{},\nabla\bu_h^{})_\Omega^{}}_{=0}
+ (p_h^{}, \dive\bu_h^{})_\Omega^{}+\langle\boldsymbol{f},\bu_h^{}\rangle_\Omega^{}\big\}\nonumber\\
&= -\lim_{h\to 0} s(p_h^{},p_h^{}) + \lim_{h\to 0} \langle\boldsymbol{f},\bu_h^{}\rangle_\Omega^{} \nonumber\\
&= 0+ \langle\boldsymbol{f},\bu\rangle_\Omega^{} \nonumber\\
&=\big (|\nabla\bu|^{r-2}\nabla\bu,\nabla\bu\big)_\Omega^{},
\end{align}
and  the result follows by using that $\bu_h^{} \rightharpoonup  \bu$ in $W^{1,r}_0(\Omega)$, 
 the fact that $W^{1,r}_0(\Omega)$ is uniformly convex, and \cite[Proposition~3.32]{Brezis}.
For $r< \frac{3d}{d+2}$ we realise that  \eqref{necessary-2} in fact states that $\lim_{h\to 0}|\bu_h^{}|_{1,r,\Omega}^{}=|\bu|_{1,r,\Omega}^{}$, 
and the strong convergence of $\bu_h^{}$ to $\bu$ in $W^{1,r}_0(\Omega)^d$ follows using once again \cite[Proposition~3.32]{Brezis}.
\end{proof}

The strong convergence of the pressure is proved next. We
 begin by noticing that, thanks to Theorem~\ref{Theo:strong-conv-u} and the continuous injection 
$W^{1,r}_0(\Omega)^d \hookrightarrow L^{2\tilde{r}}(\Omega)^d$ we have that $\bu_h^{}$ converges 
strongly to $\bu$ in $L^{2\tilde{r}}(\Omega)$. Moreover, if $r\ge \frac{3d}{d+2}$ then thanks to 
\eqref{bound-L-3}  and   \eqref{strong-jump-r}, 
$\calL(\bu_h^{},p_h^{})$ also converges  strongly to $\bu$ in $L^{2\tilde{r}}(\Omega)$.

\begin{theorem} For $r\ge \frac{3d}{d+2}$, the discrete pressure $p_h^{}$ converges to $p$ strongly in $L^{\tilde{r}}_0(\Omega)$.
\end{theorem}

\begin{proof} Let $\Pi_{H}^{}$ be the projection defined in \eqref{def-PiH}. Using the triangle inequality
we get
\begin{equation}\label{strong-p-2}
\|p-p_h^{}\|_{0,\tilde{r},\Omega}^{}\le \|p-\Pi_{H}^{}(p)\|_{0,\tilde{r},\Omega}^{}+
\|\Pi_{H}^{}(p-p_h^{})\|_{0,\tilde{r},\Omega}^{} + \|\Pi_{H}^{}(p_h^{})-p_h^{}\|_{0,\tilde{r},\Omega}^{} = \circled{1}+\circled{2}+\circled{3}.
\end{equation}
First, thanks to \eqref{conver-Pih}
\begin{equation}\label{strong-p-3}
\circled{1}\to 0\qquad\textrm{as}\; h\to 0.
\end{equation}
Moreover,  
the combined use of \eqref{bound-s} and \eqref{strong-jump-r} gives
\begin{equation}\label{strong-p-4}
\circled{3} \le\, C\,s(p_h^{},p_h^{})^{\frac{1}{2}} \to 0\qquad\textrm{as}\;h\to 0.
\end{equation}

It only remains to bound $\circled{2}$. Thanks to the inf-sup condition \eqref{inf-sup-help} there exist  $\beta_{\tilde{r}}^{}>0$ and 
$\tilde{\bv}_h^{}\in\bV_h^{}$ with $|\tilde{\bv}_h^{}|_{1,\tilde{r}',\Omega}^{}=1$ such that
\begin{align}
\beta_{\tilde{r}}^{}\, \circled{2} \, &\le\, (\Pi_{H}^{}(p-p_h^{}),\dive\tilde{\bv}_h^{})_\Omega^{} \nonumber\\
& =\, (\Pi_{H}^{}(p)-p,\dive\tilde{\bv}_h^{})_\Omega^{} +
(p-p_h^{}, \dive\tilde{\bv}_h^{})_\Omega^{} + (p_h^{}-\Pi_{H}^{}(p_h^{}),\dive\tilde{\bv}_h^{})_\Omega^{}
\nonumber\\
&= \circled{4}+\circled{5}+\circled{6}. \label{strong-p-5}
\end{align}
 H\"older's inequality gives
\begin{align}
|\circled{4}|\le \circled{1}\,\|\dive\tilde{\bv}_h^{}\|_{0,\tilde{r}',\Omega}^{}\le C\,\circled{1}\to 0, \label{strong-p-6}\\
|\circled{6}| \le \circled{3}\,\|\dive\tilde{\bv}_h^{}\|_{0,\tilde{r}',\Omega}^{}\le C\, \circled{3}\to 0,\label{strong-p-7}
\end{align}
thanks to \eqref{strong-p-3} and \eqref{strong-p-4}. It only remains to bound $\circled{5}$. Using that $(\bu_h^{},p_h^{})$ solves 
\eqref{discrete-p-laplacian-1}, \eqref{discrete-p-laplacian-2} we get
\begin{equation}\label{strong-p-8}
\circled{5} = \big( |\nabla\bu_h^{}|^{r-2}\nabla\bu_h^{}-|\nabla\bu|^{r-2}\nabla\bu,\nabla\tilde{\bv}_h^{}\big)_\Omega^{}
- \big( \calL(\bu_h^{},p_h^{})\otimes \bu_h^{}-\bu\otimes\bu,\nabla\tilde{\bv}_h^{}\big)_\Omega^{}\\
=\circled{7}+\circled{8}.
\end{equation}
Using that $\bu_h^{}$ converges to $\bu$ strongly in $W^{1,r}_0(\Omega)^d$ we get
$\circled{7} \to 0$.
Finally, $\calL(\bu_h^{},p_h^{})\otimes \bu_h^{}\to \bu\otimes\bu$ in $L^{\tilde{r}}(\Omega)^{d\times d}$ and $|\tilde{v}_h^{}|_{1,\tilde{r}',\Omega}^{}=1$ giving 
$\circled{8}\to 0$ as $h\to 0$. So, $\circled{2}\to 0$, and the result follows from \eqref{strong-p-2}.
\end{proof}

%

\section{Concluding remarks}

In this work we have extended the applicability of a low-order divergence-free stabilised finite element method to incompressible non-Newtonian fluid flow models
with power-law rheology. The method is based on using a standard continuous piecewise linear finite element approximation for the velocity and piecewise constant approximation for the pressure.  
The main results of the paper are twofold: first, the method has been shown to converge to a weak solution of the boundary-value problem in the entire
range $r>\frac{2d}{d+2}$ of the power-law index $r$ within which weak solutions to the model are known to exist. Up to now this was only possible by using finite element methods based on  pointwise divergence-free continuous piecewise polynomials constructed by taking the curl of $C^1$ piecewise polynominals (an approach that is usually avoided because of the complexity of its implementation and the excessive number of unknowns at each node, particularly in three dimensions); by using Scott--Vogelius finite elements,  which are inf-sup stable on shape-regular meshes for piecewise quartic velocity fields and higher 
\cite{GS18}; or by using Guzm\'{a}n--Neilan type pointwise divergence-free rational basis functions (see \cite{DKS13} for the convergence proof in this case). With standard mixed finite element methods, with a discretely divergence-free velocity field, the range of $r$ for which convergence was shown to hold is smaller, and is restricted to $r>\frac{2d}{d+1}$; it is not known whether such standard mixed finite element methods converge for $\frac{2d}{d+2}<r \leq \frac{2d}{d+1}$ (see \cite{DKS13}).

The second main result of this paper is the proof of strong convergence of both the velocity and the pressure. To the best of our knowledge, this is the first work where such a result has been shown for this type of stabilisation; in fact, this strong convergence result is new even for $r=2$ corresponding to the case of a Newtonian fluid. To date, not many stabilised finite element methods have been proved to be convergent
under minimal regularity hypotheses, and for those for which this was achieved the discussion was restricted to the simpler situation of a Newtonian fluid $(r=2)$. In addition, the stabilising jump terms involved the complete Cauchy stress tensor rather than the jump in the pressure alone (see, e.g., \cite{BCGG13}, where dG methods were analysed), or, in the case of continuous finite element pairs, residual-based stabilisation was used (see, \cite{BS13}).

As was noted earlier, the present work is seen a proof-of-concept paper, whose aim is to showcase the applicability of this type of stabilisation to problems that are more complex than the Navier--Stokes model, and to highlight the fact that the use of the `covert' divergence-free velocity field $\calL(\bu_h^{},p_h^{})$
in the convection term allows one to prove the convergence in the whole range of values of the power-law index $r$ for which weak solutions to the model are known to exist. As such, several questions remain  open, including the following:

\begin{itemize}
\item Assumption (A1) was introduced so as to be able to define the discretely divergence-free Lipschitz truncation in the present setting. Whether this is a necessity or one may avoid the use of Lemma~\ref{Lemma2.30} altogether, and thereby dispense with Assumption (A1), is an interesting open question;

\item The discussion contained in Remark~\ref{Rem:low-order} hints at the possibility
of applying lower-order divergence-free finite elements on appropriately refined meshes without the need of stabilisation. This would require the study of the inf-sup stability of such pairs 
in the setting of the present paper; we note in this direction the recent paper \cite{FGS}, where a Scott--Vogelius pair is used on barycentrically refined meshes.

\item Most of the results presented in this work can be extended, without major difficulties, to more sophisticated explicit constitutive laws (e.g. to Carreau--Yasuda type models). In particular, power-law models such as the ones discussed in \cite[Section~3]{DMS08} can be analysed with the techniques developed in this work;

\item Finally, the extension of the results of our work to steady and unsteady implicitly-constituted models, such as the ones considered in \cite{DKS13,ST20}, where the constitutive relation can be identified with a maximal monotone $r$-graph, is the subject of ongoing research and will be presented elsewhere.

\end{itemize}

%

\section*{Acknowledgements}

The work of GRB has been partially funded by the Leverhume Trust via the Research Fellowship No. RF-2019-510.

%

\section*{Appendix: The proof of the inf-sup condition \eqref{inf-sup-help}}

We start by introducing notation that is only used in this appendix. We define the standard finite element space 
$\PP_1^0(\calT_h^{})=\{ v_h^{}:v_h^{}|_K^{}\in\PP_1^{}(K)\;\forall\, K\in\calT_h^{}\}\cap H^1_0(\Omega)$. The set of internal facets of $\calT_H^{}$ is denoted by $\calF_H^{}$.
For each $\tilde{F}\in\calF_H^{}$ we choose one unit normal to it denoted by $\boldsymbol{n}_{\tilde{F}}^{}$. Its orientation is of no importance. Finally, for $M\in \calT_H^{}$
we define the neighbourhood $\omega_M^{}=\{ M'\in\calT_H^{}: M'\cap M\not=\emptyset\}$.

 Let $\bv\in W^{1,s'}_0(\Omega)$ and let $M\in\calT_{H}^{}$. 
Let, for any internal facet $\tilde{F}$ of $\calT_H^{}$ such that $\tilde{F}\subseteq \partial M$, $\bx$ be a node of $\calT_h^{}$ that belongs to the interior of $\tilde{F}$
(the existence of such a node is guaranteed by $\hypothesis$). Associated to $\bx$, let
$b_{\tilde{F}}^{}$ be the basis function of $\PP_1^0(\calT_h^{})$ whose value at $\bx$ is one, and zero at every other node of $\calT_h^{}$.
Then, we define the mapping
\begin{align}
\rho_M^{}:W^{1,s'}(M)&\to \bV_h^{}, \nonumber\\
\bv &\mapsto \rho_M^{}(\bv)=\sum_{\tilde{F}\subseteq\partial M} \frac{1}{(1,b_{\tilde{F}}^{})_{\tilde{F}} ^{} }\left(\bv\cdot\boldsymbol{n}_{\tilde{F}}^{},1\right)_{\tilde{F}}^{}
b_{\tilde{F}} ^{}\boldsymbol{n}_{\tilde{F}}^{}.
\end{align}
This mapping is well defined since $s'>1$ and thanks to \cite[Proposition~16.1]{EG21-I} 
the integral of the normal component of $\bv$ is finite on each ${\tilde{F}} $. 
Let now $q_{H}^{}\in\calQ_{H}^{}$
and let $\bv\in W^{1,{s}'}_0(\Omega)$ be arbitrary. Denoting by  $\boldsymbol{n}^M$ the unit normal outward to $M$ and integrating by parts we obtain
\begin{align}
\left(\dive\left(\sum_{M\in\calT_{H}^{}}\rho_M^{}(\bv)\right), q_{H}^{}\right)_\Omega &=
\sum_{M\in\calT_{H}^{}} (\rho_M^{}(\bv)\cdot\boldsymbol{n}^M,q_{H}^{})_{\partial M}^{} \nonumber\\
&= \sum_{{\tilde{F}} \in \calF_{H}^{}} \frac{1}{(1,b_{\tilde{F}}^{})_{\tilde{F}}^{} }(1, b_{\tilde{F}} ^{})_{\tilde{F}} ^{}(\bv\cdot\boldsymbol{n}_{\tilde{F}} ^{},1)_{\tilde{F}} ^{}\,\llbracket q_{H}^{}\rrbracket_{\tilde{F}}^{} \nonumber\\
&= \sum_{{\tilde{F}} \in \calF_{H}^{}} (\bv\cdot\boldsymbol{n}_{\tilde{F}} ^{},1)_{\tilde{F}} ^{}\llbracket q_{H}^{}\rrbracket_{\tilde{F}}^{} \nonumber\\
&= (\dive\bv, q_{H}^{})_\Omega^{}. \label{rhoM-div}
\end{align}
In addition, using H\"older's inequality and the local
trace inequality \eqref{discrete-trace} we obtain, for all $M\in\calT_{H}^{}$
\begin{align*}
\|\rho_M^{}(\bv)\|_{0,{s}',M}^{} &\le\, \sum_{F\subseteq \partial M} \frac{1}{(1,b_{\tilde{F}}^{})_{\tilde{F}} ^{} }\left|( \bv\cdot\boldsymbol{n}_{\tilde{F}}^{}, 1)_{\tilde{F}}^{}\right| 
\|b_{\tilde{F}}^{}\|_{0,s',M}^{} \nonumber\\
& \le\, \sum_{{\tilde{F}}\subseteq \partial M}\frac{C}{|\tilde{F}|} |M|^{\frac{1}{s'}}|{\tilde{F}}|^{\frac{1}{s}}\|\bv\|_{0,s',{\tilde{F}}}^{}\nonumber\\
&\le\, C\,\sum_{{\tilde{F}}\subseteq \partial M} h_{\tilde{F}}^{1-d+\frac{d}{s'}+\frac{d-1}{s}}\|\bv\|_{0,s',{\tilde{F}}}^{} \nonumber\\
&\le\, C\,\sum_{{\tilde{F}}\subseteq \partial M} h_{\tilde{F}}^{\frac{1}{s'}}\big( h_M^{-\frac{1}{s'}}\|\bv\|_{0,s',M}^{}
+  h_M^{1-\frac{1}{s'}} |\bv|_{1,s',M}^{}\big) \nonumber\\
&\le\, C\,\sum_{{\tilde{F}}\subseteq \partial M} \big( \|\bv\|_{0,s',M}^{}
+  h_M^{} |\bv|_{1,s',M}^{}\big),
\end{align*}
which, using the inverse inequality \eqref{inverse} on $M$ and $h_M^{}\le Ch_K^{}$ gives
\begin{equation}\label{rho-main}
|\rho_M^{}(\bv)|_{1,s',\Omega}^{}\le\, C\,\sum_{{\tilde{F}}\subseteq \partial M} \big( h_M^{-1}\|\bv\|_{0,s',M}^{}
+  |\bv|_{1,s',M}^{}\big). 
\end{equation}

Finally, we define the Fortin operator as follows:
\begin{align}
\calI:W^{1,s'}_0(\Omega) &\to \bV_h^{}, \nonumber\\
\bv &\mapsto \calI(\bv)=S_h^{}(\bv)+\sum_{M\in\calT_{H}^{}}\rho_M^{}(\bv-S_h^{}(\bv)).
\end{align}
For every $q_{H}^{}\in\calQ_{H}^{}$ and $\bv\in W^{1,s'}_0(\Omega)$, \eqref{rhoM-div} gives
\begin{equation}\label{For-1}
(\dive \calI(\bv),q_{H}^{})_\Omega^{} = (\dive\bv, q_{H}^{})_\Omega^{}.
\end{equation}
In addition, for every $\bv\in W^{1,s'}_0(\Omega)$, using \eqref{rho-main} and 
the stability and approximation properties  of $S_h^{}$ we arrive at
\begin{align}
|\calI(\bv)|_{1,s',\Omega}^{} &\le |S_h^{}(\bv)|_{1,s',\Omega}^{}+C\left\{\sum_{M\in\calT_{H}}
h_M^{-s'}\|\bv-S_h^{}(\bv)\|_{0,s',M}^{s'}+ |\bv-S_h^{}(\bv)|_{1,s',M}^{s'}\right\}^{\frac{1}{s'}}\nonumber\\
&\le C\,|\bv|_{1,s',\Omega}^{}+C\left\{\sum_{M\in\calT_{H}}|\bv|_{1,s',\omega_M^{}}^{s'}\right\}^{\frac{1}{s'}}\nonumber\\
&\le C\,|\bv|_{1,s',\Omega}^{}.\label{For-2}
\end{align}
So, $\calI$ satisfies the requirements of a Fortin operator, which proves
the inf-sup condition \eqref{inf-sup-help}.

\bibliographystyle{plain}
\bibliography{B-Suli}

\end{document}